\documentclass[a4paper,twocolumn,english]{IEEEtran}
\usepackage[T1]{fontenc}
\usepackage{verbatim}
\usepackage{float}
\usepackage{amsthm}
\usepackage{amstext}
\makeatletter

\floatstyle{ruled}
\newfloat{algorithm}{tbp}{loa}
\providecommand{\algorithmname}{Algorithm}
\floatname{algorithm}{\protect\algorithmname}

\theoremstyle{plain}
\newtheorem{thm}{\protect\theoremname}
\theoremstyle{plain}
\newtheorem{prop}[thm]{\protect\propositionname}

\usepackage{amsmath}
\usepackage{amssymb}
\usepackage[numbers,sort]{natbib}

\usepackage{graphicx}

\usepackage{subfigure}
\usepackage{pdfsync}
\usepackage{url}
\usepackage{pdfsync}
\makeatother

\usepackage{babel}
\providecommand{\propositionname}{Proposition}
\providecommand{\theoremname}{Theorem}

\begin{document}

\title{A Unified Successive Pseudo-Convex Approximation Framework}

\author{Yang Yang and Marius Pesavento%
\thanks{Y. Yang is with Intel Deutschland GmbH, Germany (email: yang1.yang@intel.com).%
}%
\thanks{M. Pesavento is with Communication Systems Group, Darmstadt University
of Technology, Germany (email: pesavento@nt.tu-darmstadt.de).%
}
\thanks{The authors acknowledge the financial support of the Seventh Framework
Programme for Research of the European Commission under grant number
ADEL-619647 and the EXPRESS project within the DFG priority program
CoSIP (DFG-SPP 1798).}}
\maketitle
\begin{abstract}
In this paper, we propose a successive pseudo-convex approximation
algorithm to efficiently compute stationary points for a large class
of possibly nonconvex optimization problems. The stationary points
are obtained by solving a sequence of successively refined approximate
problems, each of which is much easier to solve than the original
problem. To achieve convergence, the approximate problem only needs
to exhibit a weak form of convexity, namely, pseudo-convexity. We
show that the proposed framework not only includes as special cases
a number of existing methods, for example, the gradient method and
the Jacobi algorithm, but also leads to new algorithms which enjoy
easier implementation and faster convergence speed. We also propose
a novel line search method for nondifferentiable optimization problems,
which is carried out over a properly constructed differentiable function
with the benefit of a simplified implementation as compared to state-of-the-art
line search techniques that directly operate on the original nondifferentiable
objective function. The advantages of the proposed algorithm are shown,
both theoretically and numerically, by several example applications,
namely, MIMO broadcast channel capacity computation, energy efficiency
maximization in massive MIMO systems and LASSO in sparse signal recovery.\end{abstract}

\begin{IEEEkeywords}
Energy efficiency, exact line search, LASSO, massive MIMO, MIMO broadcast
channel, nonconvex optimization, nondifferentiable optimization, successive
convex approximation.
\end{IEEEkeywords}

\section{Introduction}

In this paper, we propose an iterative algorithm to solve the following
general optimization problem:
\begin{equation}
\begin{split}\underset{\mathbf{x}}{\textrm{minimize}}\quad & f(\mathbf{x})\\
\textrm{subject to}\quad & \mathbf{x}\in\mathcal{X},
\end{split}
\label{eq:original_function}
\end{equation}
where $\mathcal{X}\subseteq\mathcal{R}^{n}$ is a closed and convex
set, and $f(\mathbf{x}):\,\mathcal{R}^{n}\rightarrow\mathcal{R}$
is a proper and differentiable function with a continuous gradient.
We assume that problem (\ref{eq:original_function}) has a solution.

Problem (\ref{eq:original_function}) also includes some class of
nondifferentiable optimization problems, if the nondifferentiable
function $g(\mathbf{x})$ is convex:
\begin{equation}
\begin{split}\underset{\mathbf{x}}{\textrm{minimize}}\quad & f(\mathbf{x})+g(\mathbf{x})\\
\textrm{subject to}\quad & \mathbf{x}\in\mathcal{X},
\end{split}
\label{eq:original_function_nonsmooth}
\end{equation}
because problem (\ref{eq:original_function_nonsmooth}) can be rewritten
into a problem with the form of (\ref{eq:original_function}) by the
help of auxiliary variables:
\begin{equation}
\begin{split}\underset{\mathbf{x},y}{\textrm{minimize}}\quad & f(\mathbf{x})+y\\
\textrm{subject to}\quad & \mathbf{x}\in\mathcal{X},\, g(\mathbf{x})\leq y.
\end{split}
\label{eq:original-function-smooth}
\end{equation}

We do not assume that $f(\mathbf{x})$ is convex, so (\ref{eq:original_function})
is in general a nonconvex optimization problem. The focus of this
paper is on the development of efficient iterative algorithms for
computing the stationary points of problem (\ref{eq:original_function}).
The optimization problem (\ref{eq:original_function}) represents
general class of optimization problems with a vast number of diverse
applications. Consider for example the sum-rate maximization in the
MIMO multiple access channel (MAC) \cite{Yu2004}, the broadcast channel
(BC) \cite{Jindal2005,Yu2006} and the interference channel (IC) \cite{Ye2003,Luo2008,Kim2011,Shi2011,Scutarib,Yang2013JSAC},
where $f(\mathbf{x})$ is the sum-rate function of multiple users
(to be maximized) while the set $\mathcal{X}$ characterizes the users'
power constraints. In the context of the MIMO IC, (\ref{eq:original_function})
is a nonconvex problem and NP-hard \cite{Luo2008}. As another example,
consider portfolio optimization in which $f(\mathbf{x})$ represents
the expected return of the portfolio (to be maximized) and the set
$\mathcal{X}$ characterizes the trading constraints \cite{Yang2013b}.
Furthermore, in sparse ($l_{1}$-regularized) linear regression, $f(\mathbf{x})$
denotes the least square function and $g(\mathbf{x})$ is the sparsity
regularization function \cite{Kim2007,Boyd2010}.

Commonly used iterative algorithms belong to the class of descent
direction methods such as the conditional gradient method and the
gradient projection method for the differentiable problem (\ref{eq:original_function})
\cite{bertsekas1999nonlinear} and the proximal gradient method for
the nondifferentiable problem (\ref{eq:original_function_nonsmooth})
\cite{Parikh2014,Combettes2011}, which often suffer from slow convergence.
To speed up the convergence, the block coordinate descent (BCD) method
that uses the notion of the nonlinear best-response has been widely
studied \cite[Sec. 2.7]{bertsekas1999nonlinear}. In particular, this
method is applicable if the constraint set of (\ref{eq:original_function})
has a Cartesian product structure $\mathcal{X}=\mathcal{X}_{1}\times\ldots\times\mathcal{X}_{K}$
such that
\begin{equation}
\begin{split}\underset{\mathbf{x}=(\mathbf{x}_{k})_{k=1}^{K}}{\textrm{minimize}}\quad & f(\mathbf{x}_{1},\ldots,\mathbf{x}_{K})\\
\textrm{subject to}\quad & \mathbf{x}_{k}\in\mathcal{X}_{k},\, k=1,\ldots,K.
\end{split}
\label{eq:original-function-cartesian}
\end{equation}
The BCD method is an iterative algorithm: in each iteration, only
one variable is updated by its best-response $\mathbf{x}_{k}^{t+1}=\arg\min_{\mathbf{x}_{k}\in\mathcal{X}_{k}}f(\mathbf{x}_{1}^{t+1},\ldots,\mathbf{x}_{k-1}^{t+1},\mathbf{x}_{k},\mathbf{x}_{k+1}^{t},\ldots,\mathbf{x}_{K}^{t})$
(i.e., the point that minimizes $f(\mathbf{x})$ with respect to (w.r.t.)
the variable $\mathbf{x}_{k}$ only while the remaining variables
are fixed to their values of the preceding iteration) and the variables
are updated sequentially. This method and its variants have been successfully
adopted to many practical problems \cite{Yu2004,Kim2011,Shi2011,Razaviyayn2013,Yang2013b}.

When the number of variables is large, the convergence speed of the
BCD method may be slow due to the sequential nature of the update.
A parallel variable update based on the best-response seems attractive
as a mean to speed up the updating procedure, however, the convergence
of a parallel best-response algorithm is only guaranteed under rather
restrictive conditions, c.f. the diagonal dominance condition on the
objective function $f(\mathbf{x}_{1},\ldots,\mathbf{x}_{K})$ \cite{Bertsekas},
which is not only difficult to satisfy but also hard to verify. If
$f(\mathbf{x}_{1},\ldots,\mathbf{x}_{K})$ is convex, the parallel
algorithms converge if the stepsize is inversely proportional to the
number of block variables $K$. This choice of stepsize, however,
tends to be overly conservative in systems with a large number of
block variables and inevitably slows down the convergence \cite{Jindal2005,He2011,Yang2013b}.

A recent progress in parallel algorithms has been made in \cite{Scutarib,Yang2013JSAC,Scutari_BigData,Scutari_Nonconvex},
in which it was shown that the stationary point of (\ref{eq:original_function})
can be found by solving a sequence of successively refined \emph{approximate
problems} of the original problem (\ref{eq:original_function}), and
convergence to a stationary point is established if, among other conditions,
the approximate function (the objective function of the approximate
problem) and stepsizes are properly selected. The parallel algorithms
proposed in \cite{Scutarib,Yang2013JSAC,Scutari_BigData,Scutari_Nonconvex}
are essentially descent direction methods. A description on how to
construct the approximate problem such that the convexity of the original
problem is preserved as much as possible is also contained in \cite{Scutarib,Yang2013JSAC,Scutari_BigData,Scutari_Nonconvex}
to achieve faster convergence than standard descent directions methods
such as classical conditional gradient method and gradient projection
method.

Despite its novelty, the parallel algorithms proposed in \cite{Scutarib,Yang2013JSAC,Scutari_BigData,Scutari_Nonconvex}
suffer from two limitations. Firstly, the approximate function must
be strongly convex, and this is usually guaranteed by artificially
adding a quadratic regularization term to the original objective function
$f(\mathbf{x})$, which however may destroy the desirable characteristic
structure of the original problem that could otherwise be exploited,
e.g., to obtain computationally efficient closed-form solutions of
the approximate problems \cite{Kim2011}. Secondly, the algorithms
require the use of a decreasing stepsize. On the one hand, a slow
decay of the stepsize is preferable to make notable progress and to
achieve satisfactory convergence speed; on the other hand, theoretical
convergence is guaranteed only when the stepsize decays fast enough.
In practice, it is a difficult task on its own to find a decay rate
for the stepsize that provides a good trade-off between convergence
speed and convergence guarantee, and current practices mainly rely
on heuristics \cite{Scutari_BigData}.

The contribution of this paper consists in the development of a novel
iterative convex approximation method to solve problem (\ref{eq:original_function}).
In particular, the advantages of the proposed iterative algorithm
are the following:

1) The approximate function of the original problem (\ref{eq:original_function})
in each iteration only needs to exhibit a weak form of convexity,
namely, pseudo-convexity. The proposed iterative method not only includes
as special cases many existing methods, for example, \cite{Ye2003,Scutarib,Scutari_BigData,Yang2013JSAC},
but also opens new possibilities for constructing approximate problems
that are easier to solve. For example, in the MIMO BC sum-rate maximization
problems (Sec. \ref{sub:MIMO-Broadcast-Channel}), the new approximate
problems can be solved in closed-form. We also show by a counterexample
that the assumption on pseudo-convexity is tight in the sense that
if it is not satisfied, the algorithm may not converge.

2) The stepsizes can be determined based on the problem structure,
typically resulting in faster convergence than in cases where constant
stepsizes \cite{Jindal2005,He2011,Yang2013b} and decreasing stepsizes
\cite{Scutarib,Scutari_BigData} are used. For example, a constant
stepsize can be used when $f(\mathbf{x})$ is given as the difference
of two convex functions as in DC programming \cite{Horst1999}. When
the objective function is nondifferentiable, we propose a new exact/successive
line search method that is carried out over a properly constructed
differentiable function. Thus it is much easier to implement than
state-of-the-art techniques that operate on the original nondifferentiable
objective function directly.

In the proposed algorithm, the exact/successive line search is used
to determine the stepsize and it can be implemented in a centralized
controller, whose existence presence is justified for particular applications,
e.g., the base station in the MIMO BC, and the portfolio manager in
multi-portfolio optimization \cite{Yang2013b}. We remark that also
in applications in which centralized controller are not admitted,
however, the line search procedure does not necessarily imply an increased
signaling burden when it is implemented in a distributed manner among
different distributed processors. For example, in the LASSO problem
studied in Sec. \ref{sub:LASSO}, the stepsize based on the exact
line search can be computed in closed-form and it does not incur any
additional signaling as in predetermined stepsizes, e.g., decreasing
stepsizes and constant stepsizes. Besides, even in cases where the
line search procedure induces additional signaling, the burden is
often fully amortized by the significant increase in the convergence
rate.

The rest of the paper is organized as follows. In Sec. \ref{sec:Preliminaries}
we introduce the mathematical background. The novel iterative method
is proposed and its convergence is analyzed in Sec. \ref{sec:Proposed-method};
its connection to several existing descent direction algorithms is
presented there. In Sec. \ref{sec:Applications}, several applications
are considered: the sum rate maximization problem of MIMO BC, the
energy efficiency maximization of a massive MIMO system to illustrate
the advantage of the proposed approximate function, and the LASSO
problem to illustrate the advantage of the proposed stepsize. The
paper is finally concluded in Sec. \ref{sec:Concluding-remarks}.

\emph{Notation: }We use $x$, $\mathbf{x}$ and $\mathbf{X}$ to denote
a scalar, vector and matrix, respectively. We use $X_{jk}$ to denote
the $(j,k)$-th element of $\mathbf{X}$; $x_{k}$ is the $k$-th
element of $\mathbf{x}$ where $\mathbf{x}=(x_{k})_{k=1}^{K}$, and
$\mathbf{x}_{-k}$ denotes all elements of $\mathbf{x}$ except $x_{k}$:
$\mathbf{x}_{-k}=(x_{j})_{j=1,j\neq k}^{K}$. We denote $\mathbf{x}^{-1}$
as the element-wise inverse of $\mathbf{x}$, i.e., $(\mathbf{x}^{-1})_{k}=1/x_{k}$.
Notation $\mathbf{x}\circ\mathbf{y}$ and $\mathbf{X}\otimes\mathbf{Y}$
denotes the Hadamard product between $\mathbf{x}$ and $\mathbf{y}$,
and the Kronecker product between $\mathbf{X}$ and $\mathbf{Y}$,
respectively. The operator $[\mathbf{x}]_{\mathbf{a}}^{\mathbf{b}}$
returns the element-wise projection of $\mathbf{x}$ onto $[\mathbf{a,b}]$:
$[\mathbf{x}]_{\mathbf{a}}^{\mathbf{b}}\triangleq\max(\min(\mathbf{x},\mathbf{b}),\mathbf{a})$,
and $\left[\mathbf{x}\right]^{+}\triangleq\left[\mathbf{x}\right]_{\mathbf{0}}$.
We denote $\left\lceil x\right\rceil $ as the smallest integer that
is larger than or equal to $x$. We denote $\mathbf{d}(\mathbf{X})$
as the vector that consists of the diagonal elements of $\mathbf{X}$
and $\textrm{diag}(\mathbf{x})$ is a diagonal matrix whose diagonal
elements are as same as $\mathbf{x}$. We use $\mathbf{1}$ to denote
the vector whose elements are equal to 1.

\section{\label{sec:Preliminaries}Preliminaries on Descent Direction Method
\protect \\
and Convex Functions}

In this section, we introduce the basic definitions and concepts that
are fundamental in the development of the mathematical formalism used
in the rest of the paper.

\textbf{Stationary point}.\emph{ }A point $\mathbf{y}\in\mathcal{X}$
is a stationary point of (\ref{eq:original_function}) if
\begin{equation}
(\mathbf{x-y})^{T}\nabla f(\mathbf{y})\geq0,\;\forall\,\mathbf{x}\in\mathcal{X}.\label{eq:stationary-point}
\end{equation}
Condition (\ref{eq:stationary-point}) is the necessary condition
for local optimality of the variable $\mathbf{y}$. For nonconvex
problems, where global optimality conditions are difficult to establish,
the computation of stationary points of the optimization problem (\ref{eq:original_function})
is generally desired. If (\ref{eq:original_function}) is convex,
stationary points coincide with (globally) optimal points and condition
(\ref{eq:stationary-point}) is also sufficient for $\mathbf{y}$
to be (globally) optimal.

\textbf{Descent direction}.\emph{ }The vector $\mathbf{d}^{t}$ is
a descent direction of the function $f(\mathbf{x})$ at $\mathbf{x}=\mathbf{x}^{t}$
if
\begin{equation}
\nabla f(\mathbf{x}^{t})^{T}\mathbf{d}^{t}<0.\label{eq:def-descent-direction}
\end{equation}
If (\ref{eq:def-descent-direction}) is satisfied, the function $f(\mathbf{x})$
can be decreased when $\mathbf{x}$ is updated from $\mathbf{x}^{t}$
along direction $\mathbf{d}^{t}$. This is because in the Taylor expansion
of $f(\mathbf{x})$ around $\mathbf{x}=\mathbf{x}^{t}$ is given by:
\[
f(\mathbf{x}^{t}+\gamma\mathbf{d}^{t})=f(\mathbf{x}^{t})+\gamma\nabla f(\mathbf{x}^{t})^{T}\mathbf{d}^{t}+o(\gamma),
\]
where the first order term is negative in view of (\ref{eq:def-descent-direction}).
For sufficiently small $\gamma$, the first order term dominates all
higher order terms. More rigorously, if $\mathbf{d}^{t}$ is a descent
direction, there exists a $\bar{\gamma}^{t}>0$ such that \cite[8.2.1]{Ortega&Rheinboldt}
\[
f(\mathbf{x}^{t}+\gamma\mathbf{d}^{t})<f(\mathbf{x}^{t}),\forall\gamma\in(0,\bar{\gamma}^{t}).
\]
Note that the converse is not necessarily true, i.e., $f(\mathbf{x}^{t+1})<f(\mathbf{x}^{t})$
for arbitrary functions $f(\mathbf{x})$ does not necessarily imply
that $\mathbf{x}^{t+1}-\mathbf{x}^{t}$ is a descent direction of
$f(\mathbf{x})$ at $\mathbf{x}=\mathbf{x}^{t}$.

\textbf{Quasi-convex function}.\emph{ }A function $h(\mathbf{x})$
is quasi-convex if for any $\alpha\in[0,1]$:
\[
h((1-\alpha)\mathbf{x}+\alpha\mathbf{y})\leq\max(h(\mathbf{x}),h(\mathbf{y})),\,\forall\,\mathbf{x,y}\in\mathcal{X}.
\]
A locally optimal point $\mathbf{y}$ of a quasi-convex function $h(\mathbf{x})$
over a convex set $\mathcal{X}$ is also globally optimal, i.e.,
\[
h(\mathbf{x})\geq h(\mathbf{y}),\forall\,\mathbf{x}\in\mathcal{X}.
\]

\textbf{Pseudo-convex function}.\emph{ }A function $h(\mathbf{x})$
is pseudo-convex if \cite{Mangasarian_NonlinearProgramming}
\[
\nabla h(\mathbf{x})^{T}(\mathbf{y-x})\geq0\Longrightarrow h(\mathbf{y})\geq h(\mathbf{x}),\;\forall\mathbf{x},\mathbf{y}\in\mathcal{X}.
\]
Another equivalent definition of pseudo-convex functions is also useful
in our context \cite{Mangasarian_NonlinearProgramming}:
\begin{equation}
h(\mathbf{y})<h(\mathbf{x})\Longrightarrow\nabla h(\mathbf{x})^{T}(\mathbf{y-x})<0.\label{eq:def-pseudo-convex-function}
\end{equation}
In other words, $h(\mathbf{y})<h(\mathbf{x})$ implies that $\mathbf{y}-\mathbf{x}$
is a descent direction of $h(\mathbf{x})$. A pseudo-convex function
is also quasi-convex \cite[Th. 9.3.5]{Mangasarian_NonlinearProgramming},
and thus any locally optimal points of pseudo-convex functions are
also globally optimal.

\textbf{Convex function}.\emph{ }A function $h(\mathbf{x})$ is convex
if
\[
h(\mathbf{y})\geq h(\mathbf{x})+\nabla h(\mathbf{x})^{T}(\mathbf{y-x}),\,\forall\,\mathbf{x},\mathbf{y}\in\mathcal{X}.
\]
It is strictly convex if the above inequality is satisfied with strict
inequality whenever $\mathbf{x\neq y}$. It is easy to see that a
convex function is pseudo-convex.

\begin{figure}[t]
\includegraphics[bb=0bp 565bp 585bp 810bp,clip,scale=0.425]{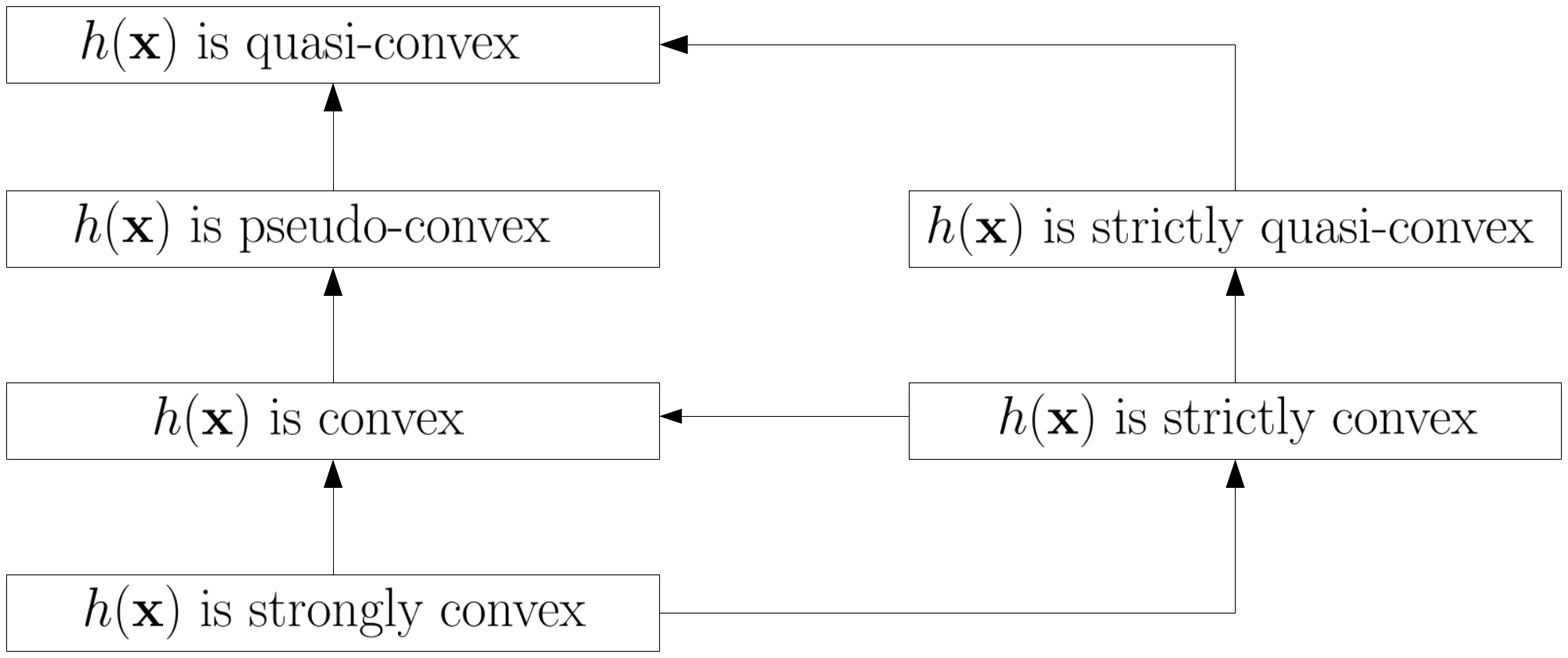}\protect\caption{\label{fig:Relationship-of-functions}Relationship of functions with
different degree of convexity}

\vspace{-1em}
\end{figure}

\textbf{Strongly convex functions}.\emph{ }A function $h(\mathbf{x})$
is strongly convex with constant $a$ if
\[
h(\mathbf{y})\geq h(\mathbf{x})+\nabla h(\mathbf{x})^{T}(\mathbf{y-x})+{\textstyle \frac{a}{2}}\left\Vert \mathbf{y}-\mathbf{x}\right\Vert _{2}^{2},\,\forall\,\mathbf{x},\mathbf{y}\in\mathcal{X},
\]
for some positive constant $a$. The relationship of functions with
different degree of convexity is summarized in Figure \ref{fig:Relationship-of-functions}
where the arrow denotes implication in the direction of the arrow.

\section{\label{sec:Proposed-method}The Proposed Successive Pseudo-Convex
Approximation Algorithm}

In this section, we propose an iterative algorithm that solves (\ref{eq:original_function})
as a sequence of successively refined approximate problems, each of
which is much easier to solve than the original problem (\ref{eq:original_function}),
e.g., the approximate problem can be decomposed into independent subproblems
that might even exhibit closed-form solutions.

In iteration $t$, let $\tilde{f}(\mathbf{x};\mathbf{x}^{t})$ be
the approximate function of $f(\mathbf{x})$ around the point $\mathbf{x}^{t}$.
Then the approximate problem is
\begin{equation}
\begin{split}\underset{\mathbf{x}}{\textrm{minimize}}\quad & \tilde{f}(\mathbf{x};\mathbf{x}^{t})\\
\textrm{subject to}\quad & \mathbf{x}\in\mathcal{X},
\end{split}
\label{eq:approximate-problem}
\end{equation}
and its optimal point and solution set is denoted as $\mathbb{B}\mathbf{x}^{t}$
and $\mathcal{S}(\mathbf{x}^{t})$, respectively:
\begin{equation}
\mathbb{B}\mathbf{x}^{t}\in\mathcal{S}(\mathbf{x}^{t})\triangleq\Bigl\{\mathbf{x}^{\star}\in\mathcal{X}:\;\tilde{f}(\mathbf{x}^{\star};\mathbf{x}^{t})=\underset{\mathbf{x}\in\mathcal{X}}{\min}\;\tilde{f}(\mathbf{x};\mathbf{x}^{t})\Bigr\}.\label{eq:mapping_definition}
\end{equation}
We assume that the approximate function $\tilde{f}(\mathbf{x};\mathbf{y})$
satisfies the following technical conditions:

\noindent (A1) The approximate function $\tilde{f}(\mathbf{x};\mathbf{y})$
is pseudo-convex in $\mathbf{x}$ for any given $\mathbf{y}\in\mathcal{X}$;

\noindent (A2) The approximate function $\tilde{f}(\mathbf{x};\mathbf{y})$
is continuously differentiable in $\mathbf{x}$ for any given $\mathbf{y}\in\mathcal{X}$
and continuous in $\mathbf{y}$ for any $\mathbf{x}\in\mathcal{X}$;

\noindent (A3) The gradient of $\tilde{f}(\mathbf{x};\mathbf{y})$
and the gradient of $f(\mathbf{x})$ are identical at $\mathbf{x}=\mathbf{y}$
for any $\mathbf{y}\in\mathcal{X}$, i.e., $\nabla_{\mathbf{x}}\tilde{f}(\mathbf{y};\mathbf{y})=\nabla_{\mathbf{x}}f(\mathbf{y})$;

Based on (\ref{eq:mapping_definition}), we define the mapping $\mathbb{B}\mathbf{x}$
that is used to generate the sequence of points in the proposed algorithm:
\begin{equation}
\mathcal{X}\ni\mathbf{x}\longmapsto\mathbb{B}\mathbf{x}\in\mathcal{X}.\label{eq:mapping}
\end{equation}
Given the mapping $\mathbb{B}\mathbf{x}$, the following properties
hold.
\begin{prop}
[Stationary point and descent direction]\label{prop:descent-property}
Provided that Assumptions (A1)-(A3) are satisfied: (i) A point $\mathbf{y}$
is a stationary point of (\ref{eq:original_function}) if and only
if $\mathbf{y}\in\mathcal{S}(\mathbf{y})$ defined in (\ref{eq:mapping_definition});
(ii) If $\mathbf{y}$ is not a stationary point of (\ref{eq:original_function}),
then $\mathbb{B}\mathbf{y}-\mathbf{y}$ is a descent direction of
$f(\mathbf{x})$:
\begin{equation}
\nabla f(\mathbf{y})^{T}(\mathbb{B}\mathbf{y}-\mathbf{y})<0.\label{eq:descent-direction}
\end{equation}
\end{prop}
\begin{IEEEproof}
See Appendix \ref{sec:Proof-of-Proposition-descent}.
\end{IEEEproof}
If $\mathbf{x}^{t}$ is not a stationary point, according to Proposition
\ref{prop:descent-property}, we define the vector update $\mathbf{x}^{t+1}$
in the $(t+1)$-th iteration as:
\begin{equation}
\mathbf{x}^{t+1}=\mathbf{x}^{t}+\gamma^{t}(\mathbb{B}\mathbf{x}^{t}-\mathbf{x}^{t}),\label{eq:variable-update}
\end{equation}
where $\gamma^{t}\in(0,1]$ is an appropriate stepsize that can be
determined by either the exact line search (also known as the minimization
rule) or the successive line search (also known as the Armijo rule).
Since $\mathbf{x}^{t}\in\mathcal{X}$, $\mathbb{B}\mathbf{x}^{t}\in\mathcal{X}$
and $\gamma^{t}\in(0,1]$, it follows from the convexity of $\mathcal{X}$
that $\mathbf{x}^{t+1}\in\mathcal{X}$ for all $t$.

\textbf{Exact line search}.\textbf{ }The stepsize is selected such
that the function $f(\mathbf{x})$ is decreased to the largest extent
along the descent direction $\mathbb{B}\mathbf{x}^{t}-\mathbf{x}^{t}$:
\begin{equation}
\gamma^{t}\in\underset{0\leq\gamma\leq1}{\arg\min}\; f(\mathbf{x}^{t}+\gamma(\mathbb{B}\mathbf{x}^{t}-\mathbf{x}^{t})).\label{eq:minimization-rule}
\end{equation}
With this stepsize rule, it is easy to see that if $\mathbf{x}^{t}$
is not a stationary point, then $f(\mathbf{x}^{t+1})<f(\mathbf{x}^{t})$.

In the special case that $f(\mathbf{x})$ in (\ref{eq:original_function})
is convex and $\gamma^{\star}$ nulls the gradient of $f(\mathbf{x}^{t}+\gamma(\mathbb{B}\mathbf{x}^{t}-\mathbf{x}^{t}))$,
i.e., $\nabla_{\gamma}f(\mathbf{x}^{t}+\gamma^{\star}(\mathbb{B}\mathbf{x}^{t}-\mathbf{x}^{t}))=0$,
then $\gamma^{t}$ in (\ref{eq:minimization-rule}) is simply the
projection of $\gamma^{\star}$ onto the interval $[0,1]$:
\[
\gamma^{t}=\left[\gamma^{\star}\right]_{0}^{1}=\begin{cases}
1, & \textrm{ if }\left.\nabla_{\gamma}f(\mathbf{x}^{t}+\gamma(\mathbb{B}\mathbf{x}^{t}-\mathbf{x}^{t}))\right|_{\gamma=1}\geq0,\\
0, & \textrm{ if }\left.\nabla_{\gamma}f(\mathbf{x}^{t}+\gamma(\mathbb{B}\mathbf{x}^{t}-\mathbf{x}^{t}))\right|_{\gamma=0}\leq0,\\
\gamma^{\star}, & \textrm{ otherwise}.
\end{cases}
\]
If $0\leq\gamma^{t}=\gamma^{\star}\leq1$, the constrained optimization
problem in (\ref{eq:minimization-rule}) is essentially unconstrained.
In some applications it is possible to compute $\gamma^{\star}$ analytically,
e.g., if $f(\mathbf{x})$ is quadratic as in the LASSO problem (Sec.
\ref{sub:LASSO}). Otherwise, for general convex functions, $\gamma^{\star}$
can be found efficiently by the bisection method as follows. Restricting
the function $f(\mathbf{x})$ to a line $\mathbf{x}^{t}+\gamma(\mathbb{B}\mathbf{x}^{t}-\mathbf{x}^{t})$,
the new function $f(\mathbf{x}^{t}+\gamma(\mathbb{B}\mathbf{x}^{t}-\mathbf{x}^{t}))$
is convex in $\gamma$ \cite{boyd2004convex}. It thus follows that
$\nabla_{\gamma}f(\mathbf{x}^{t}+\gamma(\mathbb{B}\mathbf{x}^{t}-\mathbf{x}^{t}))<0$
if $\gamma<\gamma^{\star}$ and $\nabla_{\gamma}f(\mathbf{x}^{t}+\gamma(\mathbb{B}\mathbf{x}^{t}-\mathbf{x}^{t}))>0$
if $\gamma>\gamma^{\star}$. Given an interval $[\gamma_{\textrm{low}},\gamma_{\textrm{up}}]$
containing $\gamma^{\star}$ (the initial value of $\gamma_{\textrm{low}}$
and $\gamma_{\textrm{up}}$ is 0 and 1, respectively), set $\gamma_{\textrm{mid}}=(\gamma_{\textrm{low}}+\gamma_{\textrm{up}})/2$
and refine $\gamma_{\textrm{low}}$ and $\gamma_{\textrm{up}}$ according
to the following rule:
\[
\begin{cases}
\gamma_{\textrm{low}}=\gamma_{\textrm{mid}}, & \textrm{ if }\nabla_{\gamma}f(\mathbf{x}^{t}+\gamma_{\textrm{mid}}(\mathbb{B}\mathbf{x}^{t}-\mathbf{x}^{t}))>0,\\
\gamma_{\textrm{up}}=\gamma_{\textrm{mid}}, & \textrm{ if }\nabla_{\gamma}f(\mathbf{x}^{t}+\gamma_{\textrm{mid}}(\mathbb{B}\mathbf{x}^{t}-\mathbf{x}^{t}))<0.
\end{cases}
\]
The procedure is repeated for finite times until the gap $\gamma_{\textrm{up}}-\gamma_{\textrm{low}}$
is smaller than a prescribed precision.

\textbf{Successive line search}.\textbf{ }If no structure in $f(\mathbf{x})$
(e.g., convexity) can be exploited to efficiently compute $\gamma^{t}$
according to the exact line search (\ref{eq:minimization-rule}),
the successive line search can instead be employed: given scalars
$0<\alpha<1$ and $0<\beta<1$, the stepsize $\gamma^{t}$ is set
to be $\gamma^{t}=\beta^{m_{t}}$, where $m_{t}$ is the smallest
nonnegative integer $m$ satisfying the following inequality:
\begin{equation}
f(\mathbf{x}^{t}+\beta^{m}(\mathbb{B}\mathbf{x}^{t}-\mathbf{x}^{t}))\leq f(\mathbf{x}^{t})+\alpha\beta^{m}\nabla f(\mathbf{x}^{t})^{T}(\mathbb{B}\mathbf{x}^{t}-\mathbf{x}^{t}).\label{eq:armijo-rule}
\end{equation}
Note that the existence of a finite $m_{t}$ satisfying (\ref{eq:armijo-rule})
is always guaranteed if $\mathbb{B}\mathbf{x}^{t}-\mathbf{x}^{t}$
is a descent direction at $\mathbf{x}^{t}$ and $\nabla f(\mathbf{x}^{t})^{T}(\mathbb{B}\mathbf{x}^{t}-\mathbf{x}^{t})<0$
\cite{bertsekas1999nonlinear}, i.e., from Proposition \ref{prop:descent-property}
inequality (\ref{eq:armijo-rule}) always admits a solution.

The algorithm is formally summarized in Algorithm \ref{alg:Successive-approximation-method}
and its convergence properties are given in the following theorem.

\begin{algorithm}[t]
\textbf{Data: }$t=0$ and $\mathbf{x}^{0}\in\mathcal{X}$.

Repeat the following steps until convergence:
\begin{enumerate}
\item [\textbf{S1:}] Compute $\mathbb{B}\mathbf{x}^{t}$ using (\ref{eq:mapping_definition}).
\item [\textbf{S2:}] Compute $\gamma^{t}$ by the exact line search (\ref{eq:minimization-rule})
or the successive line search (\ref{eq:armijo-rule}).
\item [\textbf{S3:}] Update $\mathbf{x}^{t+1}$ according to (\ref{eq:variable-update})
and set $t\leftarrow t+1$.
\end{enumerate}
\protect\caption{\label{alg:Successive-approximation-method}The iterative convex approximation
algorithm for differentiable problem (\ref{eq:original_function})}
\end{algorithm}

\begin{thm}
[Convergence to a stationary point]\label{thm:convergence}Consider
the sequence $\left\{ \mathbf{x}^{t}\right\} $ generated by Algorithm
\ref{alg:Successive-approximation-method}. Provided that Assumptions
(A1)-(A3) as well as the following assumptions are satisfied:

\begin{enumerate}

\item[\emph{(A4)}] The solution set $\mathcal{S}(\mathbf{x}^{t})$
is nonempty for $t=1,2,\ldots$;

\item[\emph{(A5)}] Given any convergent subsequence $\left\{ \mathbf{x}^{t}\right\} _{t\in\mathcal{T}}$
where $\mathcal{T}\subseteq\left\{ 1,2,\ldots\right\} $, the sequence
$\left\{ \mathbb{B}\mathbf{x}^{t}\right\} _{t\in\mathcal{T}}$ is
bounded.

\end{enumerate}Then any limit point of $\left\{ \mathbf{x}^{t}\right\} $
is a stationary point of (\ref{eq:original_function}).\end{thm}
\begin{IEEEproof}
See Appendix \ref{appendix:Proof-of-Theorem}.
\end{IEEEproof}
In the following we discuss some properties of the proposed Algorithm
\ref{alg:Successive-approximation-method}.

\textbf{On the conditions (A1)-(A5)}.\emph{ }The only requirement
on the convexity of the approximate function $\tilde{f}(\mathbf{x};\mathbf{x}^{t})$
is that it is pseudo-convex, cf. (A1). To the best of our knowledge,
these are the weakest conditions for descent direction methods available
in the literature. As a result, it enables the construction of new
approximate functions that can often be optimized more easily or even
in closed-form, resulting in a significant reduction of the computational
cost. Assumptions (A2)-(A3) represent standard conditions for successive
convex approximation techniques and are satisfied for many existing
approximation functions, cf. Sec. \ref{sec:Special-Cases}. Sufficient
conditions for Assumptions (A4)-(A5) are that either the feasible
set $\mathcal{X}$ in (\ref{eq:approximate-problem}) is bounded or
the approximate function in (\ref{eq:approximate-problem}) is strongly
convex \cite{Robinson1974}. We show that how these assumptions are
satisfied in popular applications considered in Sec. \ref{sec:Applications}.

\textbf{On the pseudo-convexity of the approximate function}.\emph{
}Assumption (A1) is tight in the sense that if it is not satisfied,
Proposition \ref{prop:descent-property} may not hold. Consider the
following simple example: $f(x)=x^{3}$, where $-1\leq x\leq1$ and
the point $x^{t}=0$ at iteration $t$. Choosing the approximate function
$\tilde{f}(x;x^{t})=x^{3},$ which is quasi-convex but not pseudo-convex,
all assumptions except (A1) are satisfied. It is easy to see that
$\mathbb{B}x^{t}=-1$, however $(\mathbb{B}x^{t}-x^{t})\nabla f(x^{t})=(-1-0)\cdot0=0,$
and thus $\mathbb{B}x^{t}-x^{t}$ is not a descent direction, i.e.,
inequality (\ref{eq:descent-direction}) in Proposition \ref{prop:descent-property}
is violated.

\textbf{On the stepsize. }The stepsize can be determined in a more
straightforward way if $\tilde{f}(\mathbf{x};\mathbf{x}^{t})$ is
a global upper bound of $f(\mathbf{x})$ that is exact at $\mathbf{x}=\mathbf{x}^{t}$,
i.e., assume that
\begin{enumerate}
\item [(A6)]$\tilde{f}(\mathbf{x};\mathbf{x}^{t})\geq f(\mathbf{x})$ and
$\tilde{f}(\mathbf{x}^{t};\mathbf{x}^{t})=f(\mathbf{x}^{t})$,
\end{enumerate}
then Algorithm \ref{alg:Successive-approximation-method} converges
under the choice $\gamma^{t}=1$ which results in the update $\mathbf{x}^{t+1}=\mathbb{B}\mathbf{x}^{t}$.
To see this, we first remark that $\gamma^{t}=1$ must be an optimal
point of the following problem:
\begin{equation}
1\in\underset{0\leq\gamma\leq1}{\textrm{argmin}}\;\tilde{f}(\mathbf{x}^{t}+\gamma(\mathbb{B}\mathbf{x}^{t}-\mathbf{x}^{t});\mathbf{x}^{t}),\label{eq:exact}
\end{equation}
otherwise the optimality of $\mathbb{B}\mathbf{x}^{t}$ is contradicted,
cf. (\ref{eq:mapping_definition}). At the same time, it follows from
Proposition \ref{prop:descent-property} that $\nabla\tilde{f}(\mathbf{x}^{t};\mathbf{x}^{t})^{T}(\mathbb{B}\mathbf{x}^{t}-\mathbf{x}^{t})<0$.
The successive line search over $\tilde{f}(\mathbf{x}^{t}+\gamma(\mathbb{B}\mathbf{x}^{t}-\mathbf{x}^{t}))$
thus yields a nonnegative and finite integer $m_{t}$ such that for
some $0<\alpha<1$ and $0<\beta<1$:
\begin{align}
\tilde{f}(\mathbb{B}\mathbf{x}^{t};\mathbf{x}^{t}) & \leq\tilde{f}(\mathbf{x}^{t}+\beta^{m_{t}}(\mathbb{B}\mathbf{x}^{t}-\mathbf{x}^{t});\mathbf{x}^{t})\nonumber \\
 & \leq\tilde{f}(\mathbf{x}^{t})+\alpha\beta^{m_{t}}\nabla\tilde{f}(\mathbf{x}^{t};\mathbf{x}^{t})^{T}(\mathbb{B}\mathbf{x}^{t}-\mathbf{x}^{t})\nonumber \\
 & =f(\mathbf{x}^{t})+\alpha\beta^{m_{t}}\nabla f(\mathbf{x}^{t})^{T}(\mathbb{B}\mathbf{x}^{t}-\mathbf{x}^{t}),\label{eq:successive-argument}
\end{align}
where the second inequality comes from the definition of successive
line search {[}cf. (\ref{eq:armijo-rule}){]} and the last equality
follows from Assumptions (A3) and (A6). Invoking Assumption (A6) again,
we obtain
\begin{equation}
\left.f(\mathbf{x}^{t+1})\leq f(\mathbf{x}^{t})+\alpha\beta^{m_{t}}\nabla f(\mathbf{x}^{t})^{T}(\mathbb{B}\mathbf{x}^{t}-\mathbf{x}^{t})\right|_{\mathbf{x}^{t+1}=\mathbb{B}\mathbf{x}^{t}}.\label{eq:inexact}
\end{equation}
The proof of Theorem \ref{thm:convergence} can be used verbatim to
prove the convergence of Algorithm \ref{alg:Successive-approximation-method}
with a constant stepsize $\gamma^{t}=1$.

\subsection{Nondifferentiable Optimization Problems}

In the following we show that the proposed Algorithm \ref{alg:Successive-approximation-method}
can also be applied to solve problem (\ref{eq:original-function-smooth}),
and its equivalent formulation (\ref{eq:original_function_nonsmooth})
which contains a nondifferentiable objective function. Suppose that
$\tilde{f}(\mathbf{x};\mathbf{x}^{t})$ is an approximate function
of $f(\mathbf{x})$ in (\ref{eq:original-function-smooth}) around
$\mathbf{x}^{t}$ and it satisfies Assumptions (A1)-(A3). Then the
approximation of problem (\ref{eq:original-function-smooth}) around
$(\mathbf{x}^{t},y^{t})$ is
\begin{equation}
(\mathbb{B}\mathbf{x}^{t},y^{\star}(\mathbf{x}^{t}))\triangleq\underset{(\mathbf{x},y):\mathbf{x}\in\mathcal{X},g(\mathbf{x})\leq y}{\arg\min}\tilde{f}(\mathbf{x};\mathbf{x}^{t})+y.\label{eq:nonsmooth-approximate-problem}
\end{equation}
That is, we only need to replace the differentiable function $f(\mathbf{x})$
by its approximate function $\tilde{f}(\mathbf{x};\mathbf{x}^{t})$.
To see this, it is sufficient to verify Assumption (A3) only:
\[
\begin{split}\nabla_{\mathbf{x}}(\tilde{f}(\mathbf{x}^{t};\mathbf{x}^{t})+y) & =\nabla_{\mathbf{x}}(f(\mathbf{x}^{t})+y^{t}),\\
\nabla_{y}(\tilde{f}(\mathbf{x}^{t};\mathbf{x}^{t})+y) & =\nabla_{y}(f(\mathbf{x}^{t})+y)=1.
\end{split}
\]
 Based on the exact line search, the stepsize $\gamma^{t}$ in this
case is given as
\begin{equation}
\gamma^{t}\in\underset{0\leq\gamma\leq1}{\textrm{argmin}}\bigl\{\negthinspace f(\mathbf{x}^{t}+\gamma(\mathbb{B}\mathbf{x}^{t}-\mathbf{x}^{t}))+y^{t}+\gamma(y^{\star}(\mathbf{x}^{t})-y^{t}))\negthinspace\bigr\},\label{eq:nonsmooth-stepsize-exact}
\end{equation}
where $y^{t}\geq g(\mathbf{x}^{t})$. Then the variables $\mathbf{x}^{t+1}$
and $y^{t+1}$ are defined as follows:\begin{subequations}\label{eq:nonsmooth-update-0}
\begin{align}
\mathbf{x}^{t+1} & =\mathbf{x}^{t}+\gamma^{t}(\mathbb{B}\mathbf{x}^{t}-\mathbf{x}^{t}),\label{eq:nonsmooth-update-1}\\
y^{t+1} & =y^{t}+\gamma^{t}(y^{\star}(\mathbf{x}^{t})-y^{t}).\label{eq:nonsmooth-update-2}
\end{align}
\end{subequations}The convergence of Algorithm \ref{alg:Successive-approximation-method}
with $(\mathbb{B}\mathbf{x}^{t},y^{\star}(\mathbf{x}^{t}))$ and $\gamma^{t}$
given by (\ref{eq:nonsmooth-approximate-problem})-(\ref{eq:nonsmooth-stepsize-exact})
directly follows from Theorem \ref{thm:convergence}.

The point $y^{t+1}$ given in (\ref{eq:nonsmooth-update-2}) can be
further refined:
\begin{align*}
f(\mathbf{x}^{t+1})+y^{t+1} & =f(\mathbf{x}^{t+1})+y^{t}+\gamma^{t}(y^{\star}(\mathbf{x}^{t})-y^{t})\\
 & \geq f(\mathbf{x}^{t+1})+g(\mathbf{x}^{t})+\gamma^{t}(g(\mathbb{B}\mathbf{x}^{t})-g(\mathbf{x}^{t}))\\
 & \geq f(\mathbf{x}^{t+1})+g((1-\gamma^{t})\mathbf{x}^{t}+\gamma^{t}\mathbb{B}\mathbf{x}^{t})\\
 & =f(\mathbf{x}^{t+1})+g(\mathbf{x}^{t+1}),
\end{align*}
where the first and the second inequality comes from the fact that
$y^{t}\geq g(\mathbf{x}^{t})$ as well as $y^{\star}(\mathbf{x}^{t})=g(\mathbb{B}\mathbf{x}^{t})$
and Jensen's inequality of convex functions $g(\mathbf{x})$ \cite{boyd2004convex},
respectively. Since $y^{t+1}\geq g(\mathbf{x}^{t+1})$ by definition,
the point $(\mathbf{x}^{t+1},g(\mathbf{x}^{t+1}))$ always yields
a lower value of $f(\mathbf{x})+y$ than $(\mathbf{x}^{t+1},y^{t+1})$
while $(\mathbf{x}^{t+1},g(\mathbf{x}^{t+1}))$ is still a feasible
point for problem (\ref{eq:original-function-smooth}). The update
(\ref{eq:nonsmooth-update-2}) is then replaced by the following enhanced
rule:
\begin{equation}
y^{t+1}=g(\mathbf{x}^{t+1}).\label{eq:nonsmooth-update-3}
\end{equation}
Algorithm \ref{alg:Successive-approximation-method} with $\mathbb{B}\mathbf{x}^{t}$
given in (\ref{eq:nonsmooth-update-1}) and $y^{t+1}$ given in (\ref{eq:nonsmooth-update-3})
still converges to a stationary point of (\ref{eq:original-function-smooth}).

The notation in (\ref{eq:nonsmooth-approximate-problem})-(\ref{eq:nonsmooth-stepsize-exact})
can be simplified by removing the auxiliary variable $\mathbf{y}$:
$\mathbb{B}\mathbf{x}^{t}$ in (\ref{eq:nonsmooth-approximate-problem})
can be equivalently written as
\begin{equation}
\mathbb{B}\mathbf{x}^{t}=\underset{\mathbf{x}\in\mathcal{X}}{\arg\min}\bigl\{\tilde{f}(\mathbf{x};\mathbf{x}^{t})+g(\mathbf{x})\bigr\}\label{eq:nonsmooth-approximate-problem-1}
\end{equation}

\noindent and combining (\ref{eq:nonsmooth-stepsize-exact}) and
(\ref{eq:nonsmooth-update-3}) yields
\begin{equation}
\gamma^{t}\in\underset{0\leq\gamma\leq1}{\textrm{argmin}}\bigl\{ f(\mathbf{x}^{t}+\gamma(\mathbb{B}\mathbf{x}^{t}-\mathbf{x}^{t}))+\gamma(g(\mathbb{B}\mathbf{x}^{t})-g(\mathbf{x}^{t}))\bigr\}.\label{eq:nonsmooth-stepsize-1}
\end{equation}

In the context of the successive line search, customizing the general
definition (\ref{eq:armijo-rule}) for problem (\ref{eq:original_function_nonsmooth})
yields the choice $\gamma^{t}=\beta^{m_{t}}$ with $m_{t}$ being
the smallest integer that satisfies the inequality:\vspace{-1em}

\begin{align}
f(\mathbf{x}^{t}+\beta^{m}(\mathbb{B}\mathbf{x}^{t}-\mathbf{x}^{t}))-f(\mathbf{x}^{t}) & \leq\nonumber \\
\beta^{m}\bigl(\alpha\nabla f(\mathbf{x}^{t})^{T}(\mathbb{B}\mathbf{x}^{t}-\mathbf{x}^{t})+ & (\alpha-1)(g(\mathbb{B}\mathbf{x}^{t})-g(\mathbf{x}^{t}))\bigr).\label{eq:nonsmooth-stepsize-successive-1}
\end{align}
Based on the derivations above, the proposed algorithm for the nondifferentiable
problem (\ref{eq:original_function_nonsmooth}) is formally summarized
in Algorithm \ref{alg:Successive-approximation-method-nonsmooth}.

\begin{algorithm}[t]
\textbf{Data: }$t=0$ and $\mathbf{x}^{0}\in\mathcal{X}$.

Repeat the following steps until convergence:
\begin{enumerate}
\item [\textbf{S1:}] Compute $\mathbb{B}\mathbf{x}^{t}$ using (\ref{eq:nonsmooth-approximate-problem-1}).
\item [\textbf{S2:}] Compute $\gamma^{t}$ by the exact line search (\ref{eq:nonsmooth-stepsize-1})
or the successive line search (\ref{eq:nonsmooth-stepsize-successive-1}).
\item [\textbf{S3:}] Update $\mathbf{x}^{t+1}$ according to
\[
\mathbf{x}^{t+1}=\mathbf{x}^{t}+\gamma^{t}(\mathbb{B}\mathbf{x}^{t}-\mathbf{x}^{t}).
\]
Set $t\leftarrow t+1$.
\end{enumerate}
\protect\caption{\label{alg:Successive-approximation-method-nonsmooth}The iterative
convex approximation algorithm for nondifferentiable problem (\ref{eq:original_function_nonsmooth})}
\end{algorithm}

It is much easier to calculate $\gamma^{t}$ according to (\ref{eq:nonsmooth-stepsize-1})
than in state-of-the-art techniques that directly carry out the exact
line search over the original nondifferentiable objective function
in (\ref{eq:original_function_nonsmooth}) \cite[Rule E]{Patriksson1998},
i.e.,
\[
\min_{0\leq\gamma\leq1}f(\mathbf{x}^{t}+\gamma(\mathbb{B}\mathbf{x}^{t}-\mathbf{x}^{t}))+g(\mathbf{x}^{t}+\gamma(\mathbb{B}\mathbf{x}^{t}-\mathbf{x}^{t})).
\]
This is because the objective function in (\ref{eq:nonsmooth-stepsize-1})
is differentiable in $\gamma$ while state-of-the-art techniques involve
the minimization of a nondifferentiable function. If $f(\mathbf{x})$
exhibits a specific structure such as in quadratic functions, $\gamma^{t}$
can even be calculated in closed-form. This property will be exploited
to develop fast and easily implementable algorithm for the popular
LASSO problem in Sec. \ref{sub:LASSO}.

In the proposed successive line search, the left hand side of (\ref{eq:nonsmooth-stepsize-successive-1})
depends on $f(\mathbf{x})$ while the right hand side is linear in
$\beta^{m}$. The proposed variation of the successive line search
thus involves only the evaluation of the differentiable function $f(\mathbf{x})$
and its computational complexity and signaling exchange (when implemented
in a distributed manner) is thus lower than state-of-the-art techniques
(for example \cite[Rule A']{Patriksson1998}, \cite[Equations (9)-(10)]{Tseng2009},
\cite[Remark 4]{Scutari_BigData} and \cite[Algorithm 2.1]{Byrd2013}),
in which the whole nondifferentiable function $f(\mathbf{x})+g(\mathbf{x})$
must be repeatedly evaluated (for different $m$) and compared with
a certain benchmark before $m_{t}$ is found.

\subsection{\label{sec:Special-Cases}Special Cases and New Algorithms}

In this subsection, we interpret some existing methods in the context
of Algorithm \ref{alg:Successive-approximation-method} and show that
they can be considered as special cases of the proposed algorithm.

\textbf{Conditional gradient method}:\emph{ }In this iterative algorithm
for problem (\ref{eq:original_function}), the approximate function
is given as the first-order approximation of $f(\mathbf{x})$ at $\mathbf{x}=\mathbf{x}^{t}$
\cite[Sec. 2.2.2]{bertsekas1999nonlinear}, i.e.,
\begin{equation}
\tilde{f}(\mathbf{x};\mathbf{x}^{t})=\nabla f(\mathbf{x}^{t})^{T}(\mathbf{x}-\mathbf{x}^{t}).\label{eq:conditional-gradient}
\end{equation}
Then the stepsize is selected by either the exact line search or the
successive line search.

\textbf{Gradient projection method}:\emph{ }In this iterative algorithm
for problem (\ref{eq:original_function}), $\mathbb{B}\mathbf{x}^{t}$
is given by \cite[Sec. 2.3]{bertsekas1999nonlinear}
\[
\mathbb{B}\mathbf{x}^{t}=\left[\mathbf{x}^{t}-s^{t}\nabla f(\mathbf{x}^{t})\right]_{\mathcal{X}},
\]
where $s^{t}>0$ and $\left[\mathbf{x}\right]_{\mathcal{X}}$ denotes
the projection of $\mathbf{x}$ onto $\mathcal{X}$. This is equivalent
to defining $\tilde{f}(\mathbf{x};\mathbf{x}^{t})$ in (\ref{eq:mapping_definition})
as follows:
\begin{equation}
\tilde{f}(\mathbf{x};\mathbf{x}^{t})=\nabla f(\mathbf{x}^{t})^{T}(\mathbf{x}-\mathbf{x}^{t})+\frac{1}{2s^{t}}\left\Vert \mathbf{x}-\mathbf{x}^{t}\right\Vert _{2}^{2},\label{eq:gradient-projection-approximate-function}
\end{equation}
which is the first-order approximation of $f(\mathbf{x})$ augmented
by a quadratic regularization term that is introduced to improve the
numerical stability \cite{Bertsekas}. A generalization of (\ref{eq:gradient-projection-approximate-function})
is to replace the quadratic term by $(\mathbf{x}-\mathbf{x}^{t})\mathbf{H}^{t}(\mathbf{x}-\mathbf{x}^{t})$
where $\mathbf{H}^{t}\succ\mathbf{0}$ \cite{Tseng2009}.

\textbf{Proximal gradient method}: If $f(\mathbf{x})$ is convex and
has a Lipschitz continuous gradient with a constant $L$, the proximal
gradient method for problem (\ref{eq:original_function_nonsmooth})
has the following form \cite[Sec. 4.2]{Parikh2014}:
\begin{align}
\mathbf{x}^{t+1} & =\underset{\mathbf{x}}{\arg\min}\;\Bigl\{ s^{t}g(\mathbf{x})+\frac{1}{2}\bigl\Vert\mathbf{x}-(\mathbf{x}^{t}-s^{t}\nabla f(\mathbf{x}^{t}))\bigr\Vert^{2}\Bigr\}\nonumber \\
 & =\underset{\mathbf{x}}{\arg\min}\;\Bigl\{\nabla f(\mathbf{x}^{t})(\mathbf{x}-\mathbf{x}^{t})+\frac{1}{2s^{t}}\bigl\Vert\mathbf{x}-\mathbf{x}^{t}\bigr\Vert^{2}+g(\mathbf{x})\Bigr\}.\label{eq:proximal-gradient-update}
\end{align}
where $s^{t}>0$. In the context of the proposed framework (\ref{eq:nonsmooth-approximate-problem-1}),
the update (\ref{eq:proximal-gradient-update}) is equivalent to defining
$\tilde{f}(\mathbf{x};\mathbf{x}^{t})$ as follows:
\begin{equation}
\tilde{f}(\mathbf{x};\mathbf{x}^{t})=\nabla f(\mathbf{x})^{T}(\mathbf{x}-\mathbf{x}^{t})+\frac{1}{2s^{t}}\bigl\Vert\mathbf{x}-\mathbf{x}^{t}\bigr\Vert^{2}\label{eq:proximal-gradient-approximate-function}
\end{equation}
and setting the stepsize $\gamma^{t}=1$ for all $t$. According to
Theorem \ref{thm:convergence} and the discussion following Assumption
(A6), the proposed algorithm converges under a constant unit stepsize
if $\tilde{f}(\mathbf{x};\mathbf{x}^{t})$ is a global upper bound
of $f(\mathbf{x})$, which is indeed the case when $s^{t}\leq1/L$
in view of the descent lemma \cite[Prop. A.24]{bertsekas1999nonlinear}.

\textbf{Jacobi algorithm}: In problem (\ref{eq:original_function}),
if $f(\mathbf{x})$ is convex in each $\mathbf{x}_{k}$ where $k=1,\ldots,K$
(but not necessarily jointly convex in $(\mathbf{x}_{1},\ldots,\mathbf{x}_{K})$),
the approximate function is defined as \cite{Scutarib}
\begin{equation}
\tilde{f}(\mathbf{x};\mathbf{x}^{t})={\textstyle \sum_{k=1}^{K}}\bigl(f(\mathbf{x}_{k},\mathbf{x}_{-k}^{t})+{\textstyle \frac{\tau_{k}}{2}}\left\Vert \mathbf{x}_{k}-\mathbf{x}_{k}^{t}\right\Vert _{2}^{2}\bigr),\label{eq:jacobi-approximate-function}
\end{equation}
where $\tau_{k}\geq0$ for $k=1,\ldots,K$. The $k$-th component
function $f(\mathbf{x}_{k},\mathbf{x}_{-k}^{t})+{\textstyle \frac{\tau_{k}}{2}}\left\Vert \mathbf{x}_{k}-\mathbf{x}_{k}^{t}\right\Vert _{2}^{2}$
in (\ref{eq:jacobi-approximate-function}) is obtained from the original
function $f(\mathbf{x})$ by fixing all variables except $\mathbf{x}_{k}$,
i.e., $\mathbf{x}_{-k}=\mathbf{x}_{-k}^{t}$, and further adding a
quadratic regularization term. Since $\tilde{f}(\mathbf{x};\mathbf{x}^{t})$
in (\ref{eq:jacobi-approximate-function}) is convex, Assumption (A1)
is satisfied. Based on the observations that
\[
\begin{split}\nabla_{\mathbf{x}_{k}}\tilde{f}(\mathbf{x}^{t};\mathbf{x}^{t}) & =\nabla_{\mathbf{x}_{k}}(f(\mathbf{x}_{k},\mathbf{x}_{-k}^{t})+{\textstyle \frac{\tau_{k}}{2}}\left\Vert \mathbf{x}_{k}-\mathbf{x}_{k}^{t}\right\Vert _{2}^{2})\bigr|_{\mathbf{x}_{k}=\mathbf{x}_{k}^{t}}\\
 & =\nabla_{\mathbf{x}_{k}}f(\mathbf{x}_{k},\mathbf{x}_{-k}^{t})+\tau_{k}(\mathbf{x}_{k}-\mathbf{x}_{k}^{t})\bigr|_{\mathbf{x}_{k}=\mathbf{x}_{k}^{t}}\\
 & =\nabla_{\mathbf{x}_{k}}f(\mathbf{x}^{t}),
\end{split}
\]
we conclude that Assumption (A3) is satisfied by the choice of the
approximate function in (\ref{eq:jacobi-approximate-function}). The
resulting approximate problem is given by
\begin{equation}
\begin{split}\underset{\mathbf{x}=(\mathbf{x}_{k})_{k=1}^{K}}{\textrm{minimize}}\quad & {\textstyle \sum_{k=1}^{K}}(f(\mathbf{x}_{k},\mathbf{x}_{-k}^{t})+{\textstyle \frac{\tau_{k}}{2}\left\Vert \mathbf{x}_{k}-\mathbf{x}_{k}^{t}\right\Vert ^{2})}\\
\textrm{subject to}\quad & \mathbf{x}\in\mathcal{X}.
\end{split}
\label{eq:jacobi-approximate-problem}
\end{equation}
This is commonly known as the Jacobi algorithm. The structure inside
the constraint set $\mathcal{X}$, if any, may be exploited to solve
(\ref{eq:jacobi-approximate-problem}) even more efficiently. For
example, the constraint set $\mathcal{X}$ consists of separable constraints
in the form of $\sum_{k=1}^{K}h_{k}(\mathbf{x}_{k})\leq0$ for some
convex functions $h_{k}(\mathbf{x}_{k})$. Since subproblem (\ref{eq:jacobi-approximate-problem})
is convex, primal and dual decomposition techniques can readily be
used to solve (\ref{eq:jacobi-approximate-problem}) efficiently \cite{Palomar2006a}
(such an example is studied in Sec. \ref{sub:MIMO-Broadcast-Channel}).

To guarantee the convergence, the condition proposed in \cite{Yang2013JSAC}
is that $\tau_{k}>0$ for all $k$ in (\ref{eq:jacobi-approximate-function})
unless $f(\mathbf{x})$ is strongly convex in each $\mathbf{x}_{k}$.
However, the strong convexity of $f(\mathbf{x})$ in each $\mathbf{x}_{k}$
is a strong assumption that cannot always be satisfied and the additional
quadratic regularization term that is otherwise required may destroy
the convenient structure that could otherwise be exploited, as we
will show through an example application in the MIMO BC in Sec. \ref{sub:MIMO-Broadcast-Channel}.
In the case $\tau_{k}=0$, convergence of the Jacobi algorithm (\ref{eq:jacobi-approximate-problem})
is only proved when $f(\mathbf{x})$ is jointly convex in $(\mathbf{x}_{1},\ldots,\mathbf{x}_{K})$
and the stepsize is inversely proportional to the number of variables
$K$ \cite{bertsekas1999nonlinear,Jindal2005,Yang2013b}, namely,
$\gamma^{t}=1/K$. However, the resulting convergence speed is usually
slow when $K$ is large, as we will later demonstrate numerically
in Sec. \ref{sub:MIMO-Broadcast-Channel}.

With the technical assumptions specified in Theorem \ref{thm:convergence},
the convergence of the Jacobi algorithm with the approximate problem
(\ref{eq:jacobi-approximate-problem}) and successive line search
is guaranteed even when $\tau_{k}=0$. This is because $\tilde{f}(\mathbf{x};\mathbf{x}^{t})$
in (\ref{eq:jacobi-approximate-function}) is already convex when
$\tau_{k}=0$ for all $k$ and it naturally satisfies the pseudo-convexity
assumption specified by Assumption (A1).

In the case that the constraint set $\mathcal{X}$ has a Cartesian
product structure (\ref{eq:original-function-cartesian}), the subproblem
(\ref{eq:jacobi-approximate-problem}) is naturally decomposed into
$K$ sub-problems, one for each variable, which are then solved \emph{in
parallel}. In this case, the requirement in the convexity of $f(\mathbf{x})$
in each $\mathbf{x}_{k}$ can even be relaxed to pseudo-convexity
only (although the sum function $\sum_{k=1}^{K}f(\mathbf{x}_{k},\mathbf{x}_{-k}^{t})$
is not necessarily pseudo-convex in $\mathbf{x}$ as pseudo-convexity
is not preserved under nonnegative weighted sum operator), and this
leads to the following update: $\mathbb{B}\mathbf{x}^{t}=(\mathbb{B}_{k}\mathbf{x}^{t})_{k=1}^{K}$
and

\begin{equation}
\mathbb{B}_{k}\mathbf{x}^{t}\in\underset{\mathbf{x}_{k}\in\mathcal{X}_{k}}{\arg\min}\; f(\mathbf{x}_{k},\mathbf{x}_{-k}^{t}),\, k=1,\ldots,K,\label{eq:jacobi-approximate-problem-cartesian}
\end{equation}
where $\mathbb{B}_{k}\mathbf{x}^{t}$ can be interpreted as variable
$\mathbf{x}_{k}$'s best-response to other variables $\mathbf{x}_{-k}=(\mathbf{x}_{j})_{j\neq k}$
when $\mathbf{x}_{-k}=\mathbf{x}_{-k}^{t}$. The proposed Jacobi algorithm
is formally summarized in Algorithm \ref{alg:Jacobi} and its convergence
is proved in Theorem \ref{thm:Cartesian}.

\begin{algorithm}[t]
\textbf{Data: }$t=0$ and $\mathbf{x}_{k}^{0}\in\mathcal{X}_{k}$
for all $k=1,\ldots,K$.

Repeat the following steps until convergence:
\begin{enumerate}
\item [\textbf{S1:}] For $k=1,\ldots,K$, compute $\mathbb{B}_{k}\mathbf{x}^{t}$
using (\ref{eq:jacobi-approximate-problem-cartesian}).
\item [\textbf{S2:}] Compute $\gamma^{t}$ by the exact line search (\ref{eq:minimization-rule})
or the successive line search (\ref{eq:armijo-rule}).
\item [\textbf{S3:}] Update $\mathbf{x}^{t+1}$ according to
\[
\mathbf{x}_{k}^{t+1}=\mathbf{x}_{k}^{t}+\gamma^{t}(\mathbb{B}_{k}\mathbf{x}^{t}-\mathbf{x}_{k}^{t}),\forall k=1,\ldots,K.
\]
Set $t\leftarrow t+1$.
\end{enumerate}
\protect\caption{\label{alg:Jacobi}The Jacobi algorithm for problem (\ref{eq:original-function-cartesian})}
\end{algorithm}

\begin{thm}
\label{thm:Cartesian}Consider the sequence $\{\mathbf{x}^{t}\}$
generated by Algorithm \ref{alg:Jacobi}. Provided that $f(\mathbf{x})$
is pseudo-convex in $\mathbf{x}_{k}$ for all $k=1,\ldots,K$ and
Assumptions (A4)-(A5) are satisfied. Then any limit point of the sequence
generated by Algorithm \ref{alg:Jacobi} is a stationary point of
(\ref{eq:original-function-cartesian}).\end{thm}
\begin{IEEEproof}
See Appendix \ref{sec:Appendix-of-Theorem-Cartesian}.
\end{IEEEproof}
The convergence condition specified in Theorem \ref{thm:Cartesian}
relaxes those in \cite{Scutari_BigData,Scutarib}: $f(\mathbf{x})$
only needs to be pseudo-convex in each variable $\mathbf{x}_{k}$
and no regularization term is needed (i.e., $\tau_{k}=0$). To the
best of our knowledge, this is the weakest convergence condition on
Jacobi algorithms available in the literature. We will show in Sec.
\ref{sub:Energy-Efficiency-Maximization} by an example application
of the energy efficiency maximization problem in massive MIMO systems
how the weak assumption on the approximate function's convexity proposed
in Theorem \ref{thm:convergence} can be exploited to the largest
extent. Besides this, the line search usually yields much faster convergence
than the fixed stepsize adopted in \cite{Jindal2005} even under the
same approximate problem, cf. Sec. \ref{sub:MIMO-Broadcast-Channel}.

\textbf{DC algorithm}:\emph{ }If the objective function in (\ref{eq:original_function})
is the difference of two convex functions $f_{1}(\mathbf{x})$ and
$f_{2}(\mathbf{x})$:\vspace{-0.125cm}
\[
f(\mathbf{x})=f_{1}(\mathbf{x})-f_{2}(\mathbf{x}),\vspace{-0.125cm}
\]
the following approximate function can be used:\vspace{-0.125cm}
\[
\tilde{f}(\mathbf{x};\mathbf{x}^{t})=f_{1}(\mathbf{x})-(f_{2}(\mathbf{x}^{t})+\nabla f_{2}(\mathbf{x}^{t})^{T}(\mathbf{x}-\mathbf{x}^{t})).\vspace{-0.125cm}
\]
Since $f_{2}(\mathbf{x})$ is convex and $f_{2}(\mathbf{x})\geq f_{2}(\mathbf{x}^{t})+\nabla f_{2}(\mathbf{x}^{t})^{T}(\mathbf{x}-\mathbf{x}^{t})$,
Assumption (A6) is satisfied and the constant unit stepsize can be
chosen \cite{DC_programming}.

\section{\label{sec:Applications}Example Applications}

\subsection{\label{sub:MIMO-Broadcast-Channel}MIMO Broadcast Channel Capacity
Computation}

In this subsection, we study the MIMO BC capacity computation problem
to illustrate the advantage of the proposed approximate function.

Consider a MIMO BC where the channel matrix characterizing the transmission
from the base station to user $k$ is denoted by $\mathbf{H}_{k}$,
the transmit covariance matrix of the signal from the base station
to user $k$ is denoted as $\mathbf{Q}_{k}$, and the noise at each
user $k$ is an additive independent and identically distributed Gaussian
vector with unit variance on each of its elements. Then the sum capacity
of the MIMO BC is \cite{Yu2004a}
\begin{align}
\underset{\left\{ \mathbf{Q}_{k}\right\} }{\textrm{maximize}}\quad & \log\bigl|\mathbf{I}+{\textstyle \sum_{k=1}^{K}}\mathbf{H}_{k}\mathbf{Q}_{k}\mathbf{H}_{k}^{H}\bigr|\nonumber \\
\textrm{subject to}\quad & \mathbf{Q}_{k}\succeq\mathbf{0},\; k=1,\ldots,K,{\textstyle \sum_{k=1}^{K}}\textrm{tr}(\mathbf{Q}_{k})\leq P,\label{eq:MIMO-BC}
\end{align}
where $P$ is the power budget at the base station.

Problem (\ref{eq:MIMO-BC}) is a convex problem whose solution cannot
be expressed in closed-form and can only be found iteratively. To
apply Algorithm \ref{alg:Successive-approximation-method}, we invoke
(\ref{eq:jacobi-approximate-function})-(\ref{eq:jacobi-approximate-problem})
and the approximate problem at the $t$-th iteration is
\begin{align}
\underset{\left\{ \mathbf{Q}_{k}\right\} }{\textrm{maximize}}\quad & {\textstyle \sum_{k=1}^{K}}\log\bigl|\mathbf{R}_{k}(\mathbf{Q}_{-k}^{t})+\mathbf{H}_{k}\mathbf{Q}_{k}\mathbf{H}_{k}^{H}\bigr|\nonumber \\
\textrm{subject to}\quad & \mathbf{Q}_{k}\succeq\mathbf{0},\; k=1,\ldots,K,\;{\textstyle \sum_{k=1}^{K}}\textrm{tr}(\mathbf{Q}_{k})\leq P,\label{eq:MIMO-BC-approx}
\end{align}
where $\mathbf{R}_{k}(\mathbf{Q}_{-k}^{t})\triangleq\mathbf{I}+\sum_{j\neq k}\mathbf{H}_{j}\mathbf{Q}_{j}^{t}\mathbf{H}_{j}^{H}$.
The approximate function is concave in $\mathbf{Q}$ and differentiable
in both $\mathbf{Q}$ and $\mathbf{Q}^{t}$, and thus Assumptions
(A1)-(A3) are satisfied. Since the constraint set in (\ref{eq:MIMO-BC-approx})
is compact, the approximate problem (\ref{eq:MIMO-BC-approx}) has
a solution and Assumptions (A4)-(A5) are satisfied.

Problem (\ref{eq:MIMO-BC-approx}) is convex and the sum-power constraint
coupling $\mathbf{Q}_{1},\ldots,\mathbf{Q}_{K}$ is separable, so
dual decomposition techniques can be used \cite{Palomar2006a}. In
particular, the constraint set has a nonempty interior, so strong
duality holds and (\ref{eq:MIMO-BC-approx}) can be solved from the
dual domain by relaxing the sum-power constraint into the Lagrangian
\cite{boyd2004convex}:
\begin{equation}
\mathbb{B}\mathbf{Q}^{t}=\underset{(\mathbf{Q}_{k}\succeq\mathbf{0})_{k=1}^{K}}{\arg\max}\left\{ \negthickspace\negmedspace\begin{array}{l}
{\textstyle \sum_{k=1}^{K}}\log\left|\mathbf{R}_{k}(\mathbf{Q}_{-k}^{t})+\mathbf{H}_{k}\mathbf{Q}_{k}\mathbf{H}_{k}^{H}\right|\smallskip\\
-\lambda^{\star}({\textstyle \sum_{k=1}^{K}}\textrm{tr}(\mathbf{Q}_{k})-P)
\end{array}\negthickspace\negmedspace\right\} .\label{eq:MIMO-BC-dual}
\end{equation}
where $\mathbb{B}\mathbf{Q}^{t}=(\mathbb{B}_{k}\mathbf{Q}^{t})_{k=1}^{K}$
and $\lambda^{\star}$ is the optimal Lagrange multiplier that satisfies
the following conditions: $\lambda^{\star}\geq0$, $\sum_{k=1}^{K}\textrm{tr}(\mathbb{B}_{k}\mathbf{Q}^{t})-P\leq0$,
$\lambda^{\star}(\sum_{k=1}^{K}\textrm{tr}(\mathbb{B}_{k}\mathbf{Q}^{t})-P)=0$,
and can be found efficiently using the bisection method.

The problem in (\ref{eq:MIMO-BC-dual}) is uncoupled among different
variables $\mathbf{Q}_{k}$ in both the objective function and the
constraint set, so it can be decomposed into a set of smaller subproblems
which are solved in parallel: $\mathbb{B}\mathbf{Q}^{t}=(\mathbb{B}_{k}\mathbf{Q}^{t})_{k=1}^{K}$
and
\begin{equation}
\mathbb{B}_{k}\mathbf{Q}^{t}=\underset{\mathbf{Q}_{k}\succeq\mathbf{0}}{\arg\max}\bigl\{\log\left|\mathbf{R}_{k}(\mathbf{Q}_{-k}^{t})+\mathbf{H}_{k}\mathbf{Q}_{k}\mathbf{H}_{k}^{H}\right|-\lambda^{\star}\textrm{tr}(\mathbf{Q}_{k})\bigr\},\label{eq:MIMO-BC-dual-decomp}
\end{equation}
and $\mathbb{B}_{k}\mathbf{Q}^{t}$ exhibits a closed-form expression
based on the waterfilling solution \cite{Jindal2005}. Thus problem
(\ref{eq:MIMO-BC-approx}) also has a closed-form solution up to a
Lagrange multiplier that can be found efficiently using the bisection
method. With the update direction $\mathbb{B}\mathbf{Q}^{t}-\mathbf{Q}^{t}$,
the base station can implement the exact line search to determine
the stepsize using the bisection method described after (\ref{eq:minimization-rule})
in Sec. \ref{sec:Proposed-method}.

We remark that when the channel matrices $\mathbf{H}_{k}$ are rank
deficient, problem (\ref{eq:MIMO-BC-approx}) is convex but not strongly
convex, but the proposed algorithm with the approximate problem (\ref{eq:MIMO-BC-approx})
still converges. However, if the approximate function in \cite{Scutarib}
is used {[}cf. (\ref{eq:jacobi-approximate-function}){]}, an additional
quadratic regularization term must be included into (\ref{eq:MIMO-BC-approx})
(and thus (\ref{eq:MIMO-BC-dual-decomp})) to make the approximate
problem strongly convex, but the resulting approximate problem no
longer exhibits a closed-form solution and thus are much more difficult
to solve.

\textbf{Simulations. }The parameters are set as follows. The number
of users is $K=20$ and $K=100$, the number of transmit and receive
antenna is (5,4), and $P=10$ dB. The simulation results are averaged
over 20 instances.

\begin{figure}[t]
\center

\includegraphics[scale=0.67]{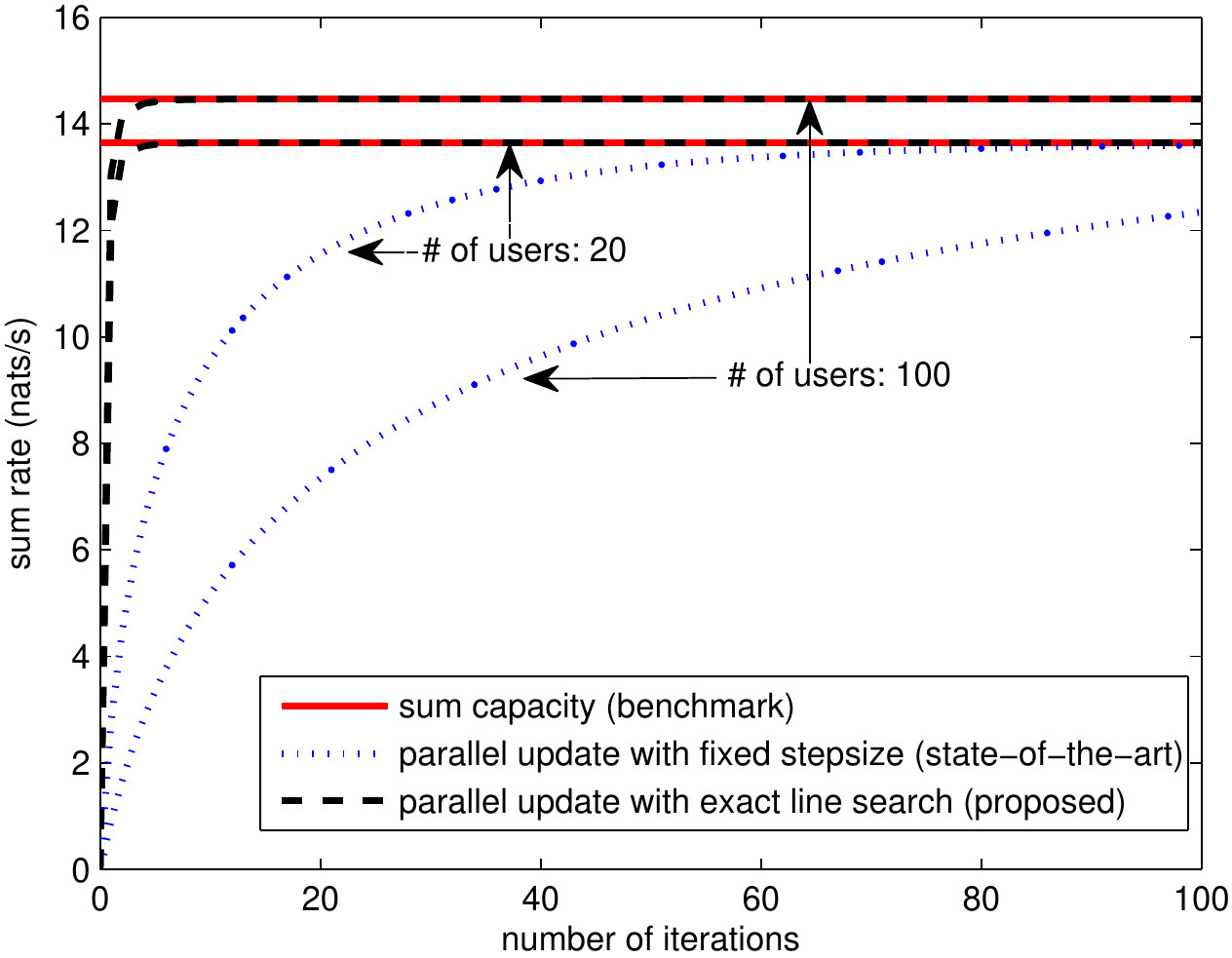}

\protect\caption{\label{fig:MIMO-BC-Sum-Capacity}MIMO BC: sum-rate versus the number
of iterations.}

\vspace{-1em}
\end{figure}

We apply Algorithm \ref{alg:Successive-approximation-method} with
approximate problem (\ref{eq:MIMO-BC-approx}) and stepsize based
on the exact line search, and compare it with the iterative algorithm
proposed in \cite{Jindal2005,He2011}, which uses the same approximate
problem (\ref{eq:MIMO-BC-approx}) but with a fixed stepsize $\gamma^{t}=1/K$
($K$ is the number of users). It is easy to see from Figure \ref{fig:MIMO-BC-Sum-Capacity}
that the proposed method converges very fast (in less than 10 iterations)
to the sum capacity, while the method of \cite{Jindal2005} requires
many more iterations. This is due to the benefit of the exact line
search applied in our algorithm over the fixed stepsize which tends
to be overly conservative. Employing the exact line search adds complexity
as compared to the simple choice of a fixed stepsize, however, since
the objective function of (\ref{eq:MIMO-BC}) is concave, the exact
line search consists in maximizing a differentiable concave function
with a scalar variable, and it can be solved efficiently by the bisection
method with affordable cost. More specifically, it takes 0.0023 seconds
to solve problem (\ref{eq:MIMO-BC-approx}) and 0.0018 seconds to
perform the exact line search (the software/hardware environment is
further specified in Sec. \ref{sub:LASSO}). Therefore, the overall
CPU time (time per iteration$\times$number of iterations) is still
dramatically decreased due to the notable reduction in the number
of iterations. Besides, in contrast to the method of \cite{Jindal2005},
increasing the number of users $K$ does not slow down the convergence,
so the proposed algorithm is scalable in large networks.

\begin{figure}[t]
\center

\includegraphics[scale=0.67]{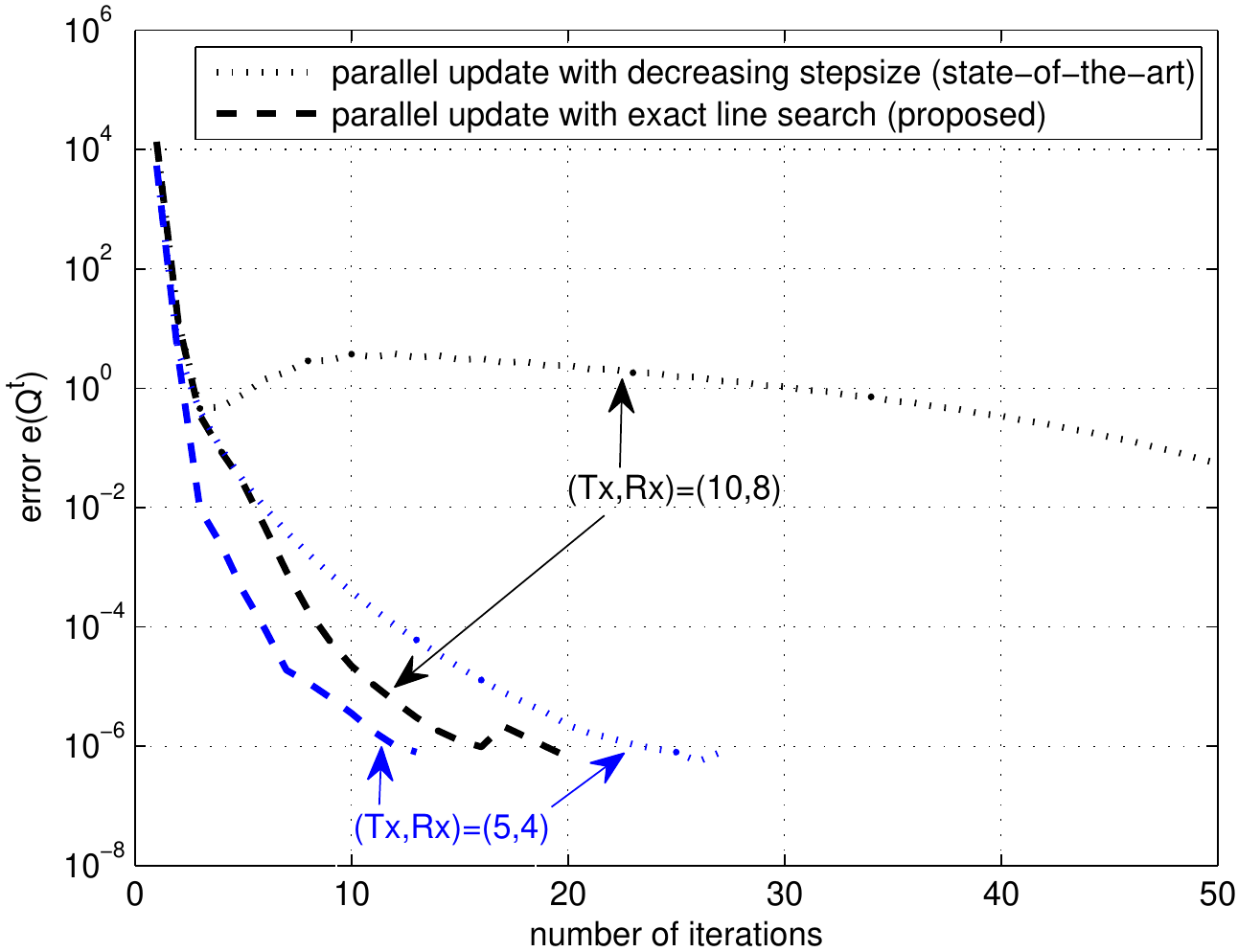}\protect\caption{\label{fig:MIMO-BC-variable-error}MIMO BC: error $e(\mathbf{Q}^{t})=\Re(\textrm{tr}(\nabla f(\mathbf{Q}^{t})(\mathbb{B}\mathbf{Q}^{t}-\mathbf{Q}^{t})))$
versus the number of iterations.}
\end{figure}

We also compare the proposed algorithm with the iterative algorithm
of \cite{Scutari_Nonconvex}, which uses the approximate problem (\ref{eq:MIMO-BC-approx})
but with an additional quadratic regularization term, cf. (\ref{eq:jacobi-approximate-function}),
where $\tau_{k}=10^{-5}$ for all $k$, and decreasing stepsizes $\gamma^{t+1}=\gamma^{t}(1-d\gamma^{t})$
where $d=0.01$ is the so-called decreasing rate that controls the
rate of decrease in the stepsize. We can see from Figure \ref{fig:MIMO-BC-variable-error}
that the convergence behavior of \cite{Scutari_Nonconvex} is rather
sensitive to the decreasing rate $d$. The choice $d=0.01$ performs
well when the number of transmit and receive antennas is 5 and 4,
respectively, but it is no longer a good choice when the number of
transmit and receive antenna increases to 10 and 8, respectively.
A good decreasing rate $d$ is usually dependent on the problem parameters
and no general rule performs equally well for all choices of parameters.

We remark once again that the complexity of each iteration of the
proposed algorithm is very low because of the existence of a closed-form
solution to the approximate problem (\ref{eq:MIMO-BC-approx}), while
the approximate problem proposed in \cite{Scutari_Nonconvex} does
not exhibit a closed-form solution and can only be solved iteratively.
Specifically, it takes $\texttt{CVX}$ (version 2.0 \cite{grant2011})
21.1785 seconds (based on the dual approach (\ref{eq:MIMO-BC-dual-decomp})
where $\lambda^{\star}$ is found by bisection). Therefore, the overall
complexity per iteration of the proposed algorithm is much lower than
that of \cite{Scutari_Nonconvex}.

\subsection{\label{sub:Energy-Efficiency-Maximization}Energy Efficiency Maximization
in Massive MIMO Systems}

In this subsection, we study the energy efficiency maximization problem
in massive MIMO systems to illustrate the advantage of the relaxed
convexity requirement of the approximate function $\tilde{f}(\mathbf{x};\mathbf{x}^{t})$
in the proposed iterative optimization approach: according to Assumption
(A1), $f(\mathbf{x};\mathbf{x}^{t})$ only needs to exhibit the pseudo-convexity
property rather than the convexity or strong convexity property that
is conventionally required.

Consider the massive MIMO network with $K$ cells and each cell serves
one user. The achievable transmission rate for user $k$ in the uplink
can be formulated into the following general form%
\footnote{We assume a single resource block. However, the extension to multiple
resource blocks is straightforward.%
}:
\begin{equation}
r_{k}(\mathbf{p})\triangleq\log\left(1+\frac{w_{kk}p_{k}}{\sigma_{k}^{2}+\phi_{k}p_{k}+\sum_{j\neq k}w_{kj}p_{j}}\right),\label{eq:EE-rate}
\end{equation}
where $p_{k}$ is the transmission power for user $k$, $\sigma_{k}^{2}$
is the covariance of the additive noise at the receiver of user $k$,
while $\phi_{k}$ and $\{w_{kj}\}_{k,j}$ are positive constants that
depend on the channel conditions only. In particular, $\phi_{k}p_{k}$
accounts for the hardware impairments, and $\sum_{j\neq k}w_{kj}p_{j}$
accounts for the interference from other users \cite{Zappone2015}.

In 5G wireless communication networks, the energy efficiency is a
key performance indicator. To address this issue, we look for the
optimal power allocation that maximizes the energy efficiency:
\begin{align}
\underset{\mathbf{p}}{\textrm{maximize}}\quad & \frac{\sum_{k=1}^{K}r_{k}(\mathbf{p})}{P_{c}+\sum_{k=1}^{K}p_{k}}\nonumber \\
\textrm{subject to}\quad & \underline{\mathbf{p}}\leq\mathbf{p}\leq\bar{\mathbf{p}},\label{eq:EE-problem-formulation}
\end{align}
where $P_{c}$ is a positive constant representing the total circuit
power dissipated in the network, $\underline{\mathbf{p}}=(\underline{p}_{k})_{k=1}^{K}$
and $\overline{\mathbf{p}}=(\overline{p}_{k})_{k=1}^{K}$ specifies
the lower and upper bound constraint, respectively.

Problem (\ref{eq:EE-problem-formulation}) is nonconvex and it is
a NP-hard problem to find a globally optimal point \cite{Luo2008}.
Therefore we aim at finding a stationary point of (\ref{eq:EE-problem-formulation})
using the proposed algorithm. To begin with, we propose the following
approximate function at $\mathbf{p}=\mathbf{p}^{t}$:
\begin{equation}
\tilde{f}(\mathbf{p};\mathbf{p}^{t})=\sum_{k=1}^{K}\frac{\tilde{r}_{k}(p_{k};\mathbf{p}^{t})}{P_{c}+p_{k}+\sum_{j\neq k}p_{j}^{t}},\label{eq:EE-approximate-function}
\end{equation}
where
\begin{equation}
\tilde{r}_{k}(p_{k};\mathbf{p}^{t})\triangleq r_{k}(p_{k},\mathbf{p}_{-k}^{t})+\sum_{j\neq k}(r_{j}(\mathbf{p}^{t})+(p_{k}-p_{k}^{t})\nabla_{p_{k}}r_{j}(\mathbf{p}^{t})).\label{eq:EE-approximate-function-individual}
\end{equation}
and
\begin{align*}
\nabla_{p_{k}}r_{j}(\mathbf{p}^{t}))= & \:\frac{w_{jk}}{\sigma_{j}^{2}+\phi_{j}p_{j}+\sum_{l=1}^{K}w_{jl}p_{l}}\\
 & -\frac{w_{jk}}{\sigma_{j}^{2}+\phi_{j}p_{j}+\sum_{l\neq j}w_{jl}p_{l}}
\end{align*}
Note that the approximate function $\tilde{f}(\mathbf{p};\mathbf{p}^{t})$
consists of $K$ component functions, one for each variable $p_{k}$,
and the $k$-th component function is constructed as follows: since
$r_{k}(\mathbf{p})$ is concave in $p_{k}$ (shown in the right column
of this page) but $r_{j}(\mathbf{p})$ is not concave in $p_{k}$
(as a matter of fact, it is convex in $p_{k}$), the concave function
$r_{k}(p_{k},\mathbf{p}_{-k}^{t})$ is preserved in $\tilde{r}_{k}(p_{k};\mathbf{p}^{t})$
in (\ref{eq:EE-approximate-function-individual}) with $\mathbf{p}_{-k}$
fixed to be $\mathbf{p}_{-k}^{t}$ while the nonconcave functions
$\{r_{j}(\mathbf{p})\}$ are linearized w.r.t. $p_{k}$ at $\mathbf{p}=\mathbf{p}^{t}$.
In this way, the partial concavity in the nonconcave function $\sum_{j=1}^{K}r_{j}(\mathbf{p})$
is preserved in $\tilde{f}(\mathbf{p};\mathbf{p}^{t})$. Similarly,
since $P_{c}+\sum_{j=1}^{K}$ in the denominator is linear in $p_{k}$,
we only set $\mathbf{p}_{-k}=\mathbf{p}_{-k}^{t}$. Note that the
division operator in the original problem (\ref{eq:EE-problem-formulation})
is kept in the approximate function $\tilde{f}(\mathbf{p};\mathbf{p}^{t})$
(\ref{eq:EE-approximate-function}). Although it will destroy the
concavity (recall that a concave function divided by a linear function
is no longer a concave function), the pseudo-concavity of $\tilde{r}_{k}(p_{k};\mathbf{p}^{t})/(P_{c}+p_{k}+\sum_{j\neq k}p_{j}^{t})$
is still preserved, as we show in two steps.

Step 1: The function $r_{k}(p_{k},\mathbf{p}_{-k}^{t})$ is concave
in $p_{k}$. For the simplicity of notation, we define two constants
$c_{1}\triangleq w_{kk}/\phi_{k}>0$ and $c_{2}\triangleq(\sigma_{k}^{2}+\sum_{j\neq k}w_{kj}p_{j}^{(t)})/\phi_{k}>0$.
The first-order derivative and second-order derivative of $r_{k}(p_{k},\mathbf{p}_{-k}^{t})$
w.r.t. $p_{k}$ are
\begin{align*}
\nabla_{p_{k}}r_{k}(p_{k},\mathbf{p}_{-k}^{t}) & =\frac{1+c_{1}}{(1+c_{1})p_{k}+c_{2}}-\frac{1}{p_{k}+c_{2}},\\
\nabla_{p_{k}}^{2}r_{k}(p_{k},\mathbf{p}_{-k}^{t}) & =-\frac{(1+c_{1})^{2}}{((1+c_{1})p_{k}+c_{2})^{2}}+\frac{1}{(p_{k}+c_{2})^{2}}\\
 & =\frac{-2c_{1}c_{2}p_{k}(1+c_{1})-(c_{1}^{2}+2c_{1})c_{2}^{2}}{((1+c_{1})p_{k}+c_{2})^{2}(p_{k}+c_{2})^{2}}.
\end{align*}
Since $\nabla_{p_{k}}^{2}r_{k}(p_{k},\mathbf{p}_{-k}^{t})<0$ when
$p\geq0$, $r_{k}(p_{k},\mathbf{p}_{-k}^{t})$ is a concave function
of $p_{k}$ in the nonnegative orthant $p_{k}\geq0$ \cite{boyd2004convex}.

Step 2: Given the concavity of $r_{k}(p_{k},\mathbf{p}_{-k}^{t})$,
the function$\tilde{r}_{k}(p_{k};\mathbf{p}^{t})$ is concave in $p_{k}$.
Since the denominator function of $\tilde{f}_{k}(p_{k};\mathbf{p}^{t})$
is a linear (and thus convex) function of $p_{k}$, it follows from
\cite[Lemma 3.8]{Olsson2007} that $\tilde{r}_{k}(p_{k};\mathbf{p}^{t})/(P_{c}+p_{k}+\sum_{j\neq k}p_{j}^{t})$
is pseudo-concave.

 \newcounter{MYtempeqncnt} \begin{figure*}[t] \normalsize \setcounter{MYtempeqncnt}{\value{equation}} \setcounter{equation}{42} \vspace*{4pt}
\begin{equation}
p_{k}(\lambda_{k}^{t,\tau})=\left[\texttt{int}_{k}(\mathbf{p}^{t})\frac{\sqrt{(2\phi_{k}+w_{kk})^{2}-4\phi_{k}\Bigl(\frac{w_{kk}}{(\pi_{k}(\mathbf{p}^{t})-\lambda_{k}^{t,\tau})\texttt{int}_{k}(\mathbf{p}^{t})}+1\Bigr)}-1}{2\phi_{k}(\pi_{k}(\mathbf{p}^{t})-\lambda_{k}^{t,\tau})(\phi_{k}+w_{kk})}\right]_{\underline{p}_{k}}^{\overline{p}_{k}},\label{eq:EE-closed-form-3}
\end{equation}
\setcounter{equation}{\value{MYtempeqncnt}} \hrulefill  \end{figure*}

Then we verify that the gradient of the approximate function and that
of the original objective function are identical at $\mathbf{p}=\mathbf{p}^{t}$.
It follows that
\begin{align*}
\left.\nabla_{p_{k}}\tilde{f}(\mathbf{p};\mathbf{p}^{t})\right|_{\mathbf{p}=\mathbf{p}^{t}} & =\left.\nabla_{p_{k}}\biggl(\frac{\tilde{r}_{k}(p_{k};\mathbf{p}^{t})}{P_{c}+p_{k}+\sum_{j\neq k}p_{j}^{t}}\biggr)\right|_{p_{k}=p_{k}^{t}}\\
 & =\sum_{j=1}^{K}\frac{\nabla_{p_{k}}r_{j}(\mathbf{p}^{t})(p_{c}+\mathbf{1}^{T}\mathbf{p}^{t})-r_{j}(\mathbf{p}^{t})}{(P_{c}+\mathbf{1}^{T}\mathbf{p}^{t})^{2}}\\
 & =\left.\nabla_{p_{k}}\Biggl(\frac{\sum_{j=1}^{K}r_{j}(\mathbf{p})}{P_{c}+\sum_{j=1}^{K}p_{j}}\Biggr)\right|_{\mathbf{p}=\mathbf{p}^{t}},\:\forall\, k.
\end{align*}
Therefore Assumption (A3) in Theorem \ref{thm:Cartesian} is satisfied.
Assumption (A2) is also satisfied because both $\tilde{r}_{k}(p_{k};\mathbf{p}^{t})$
and $p_{k}+P_{c}+\sum_{j\neq k}p_{j}^{t}$ are continuously differentiable
for any $\mathbf{p}^{t}\geq\mathbf{0}$.

Given the approximate function (\ref{eq:EE-approximate-function}),
the approximate problem in iteration $t$ is thus
\begin{equation}
\mathbb{B}\mathbf{p}^{t}=\underset{\underline{\mathbf{p}}\leq\mathbf{p}\leq\bar{\mathbf{p}}}{\arg\max}\;\sum_{k=1}^{K}\frac{\tilde{r}_{k}(p_{k};\mathbf{p}^{t})}{P_{c}+p_{k}+\sum_{j\neq k}p_{j}^{t}}.\label{eq:EE-Approximate-problem-overall}
\end{equation}
Assumptions (A4) and (A5) can be proved to hold in a similar procedure
shown in the previous example application. Since the objective function
in (\ref{eq:EE-problem-formulation}) is nonconcave, it may not be
computationally affordable to perform the exact line search. Instead,
the successive line search can be applied to calculate the stepsize.
As a result, the convergence of the proposed algorithm with approximate
problem (\ref{eq:EE-Approximate-problem-overall}) and successive
line search follows from the same line of analysis used in the proof
of Theorem \ref{thm:Cartesian}.

The optimization problem in (\ref{eq:EE-Approximate-problem-overall})
can be decomposed into independent subproblems (\ref{eq:EE-Approximate-problem-decomposition})
that can be solved in parallel:
\begin{equation}
\mathbb{B}_{k}\mathbf{p}^{t}=\underset{\underline{p}_{k}\leq p_{k}\leq\bar{p}_{k}}{\arg\max}\;\frac{\tilde{r}_{k}(p_{k};\mathbf{p}^{t})}{P_{c}+p_{k}+\sum_{j\neq k}p_{j}^{t}},\; k=1,\ldots,K,\label{eq:EE-Approximate-problem-decomposition}
\end{equation}
where $\mathbb{B}\mathbf{p}^{t}=(\mathbb{B}_{k}\mathbf{p}^{t})_{k=1}^{K}$.
As we have just shown, the numerator function and the denominator
function in (\ref{eq:EE-Approximate-problem-decomposition}) is concave
and linear, respectively, so the optimization problem in (\ref{eq:EE-Approximate-problem-decomposition})
is a fractional programming problem and can be solved by the Dinkelbach's
algorithm \cite[Algorithm 5]{Zappone2015}: given $\lambda_{k}^{t,\tau}$
($\lambda_{k}^{t,0}$ can be set to 0), the following optimization
problem in iteration $\tau+1$ is solved:\begin{subequations}\label{eq:EE-closed-form}
\begin{equation}
p_{k}(\lambda_{k}^{t,\tau})\triangleq\underset{\underline{p}_{k}\leq p_{k}\leq\overline{p}_{k}}{\arg\max}\,\tilde{r}_{k}(p_{k};\mathbf{p}^{t})-\lambda_{k}^{t,\tau}(P_{c}+p_{k}+{\textstyle \sum_{j\neq k}}p_{j}^{t}),\label{eq:EE-closed-form-1}
\end{equation}
where $\tilde{r}_{k}(p_{k};\mathbf{p}^{t})$ is the numerator function
in (\ref{eq:EE-Approximate-problem-decomposition}). The variable
$\lambda_{k}^{t,\tau}$ is then updated in iteration $\tau+1$ as
\begin{equation}
\lambda_{k}^{t,\tau+1}=\frac{\tilde{r}_{k}(p_{k}(\lambda_{k}^{t,\tau});\mathbf{p}^{t})}{P_{c}+p_{k}(\lambda_{k}^{t,\tau})+\sum_{j\neq k}p_{j}^{t}}.\label{eq:EE-closed-form-2}
\end{equation}
\end{subequations}It follows from the convergence properties of the
Dinkelbach's algorithm that
\[
\lim_{\tau\rightarrow\infty}p_{k}(\lambda_{k}^{t,\tau})=\mathbb{B}_{k}\mathbf{p}^{t}
\]
at a superlinear convergence rate. Note that problem (\ref{eq:EE-closed-form-1})
can be solved in closed-form, as $p_{k}(\lambda_{k}^{t,\tau})$ is
simply the projection of the point that sets the gradient of the objective
function in (\ref{eq:EE-closed-form-1}) to zero onto the interval
$[\underline{p}_{k},\overline{p}_{k}]$. It can be verified that finding
that point is equivalent to finding the root of a polynomial with
order 2 and it thus admits a closed-form expression. We omit the detailed
derivations and directly give the expression of $p_{k}(\lambda_{k}^{t,\tau})$
in (\ref{eq:EE-closed-form-3}) at the top of this page, where $\pi_{k}(\mathbf{p}^{t})\triangleq\sum_{j\neq k}\nabla_{p_{k}}r_{j}(\mathbf{p}^{t})$
and $\texttt{int}_{k}(\mathbf{p}^{t})\triangleq\sigma_{k}^{2}+\sum_{j\neq k}w_{kj}p_{j}^{t}$.

We finally remark that the approximate function in (\ref{eq:EE-approximate-function})
is constructed in the same spirit as \cite{Scutarib,Yang2013JSAC,Yang_stochastic}
by keeping as much concavity as possible, namely, $r_{k}(\mathbf{p})$
in the numerator and $P_{c}+\sum_{j=1}^{K}p_{j}$ in the denominator,
and linearizing the nonconcave functions only, namely, $\sum_{j\neq k}r_{j}(\mathbf{p})$
in the numerator. Besides this, the denominator function is also kept.
Therefore, the proposed algorithm is of a best-response nature and
expected to converge faster than gradient based algorithms which linearizes
the objective function $\sum_{j=1}^{K}r_{j}(\mathbf{p})/(P_{c}+\sum_{j=1}^{K}p_{j})$
in (\ref{eq:EE-problem-formulation}) completely. However, the convergence
of the proposed algorithm with the approximate problem given in (\ref{eq:EE-Approximate-problem-overall})
cannot be derived from existing works, since the approximate function
presents only a weak form of convexity, namely, the pseudo-convexity,
which is much weaker than those required in state-of-the-art convergence
analysis, e.g., uniform strong convexity in \cite{Scutarib}.

\begin{figure*}[t]
\subfigure[Number of users: $K=10$]{\includegraphics[scale=0.67]{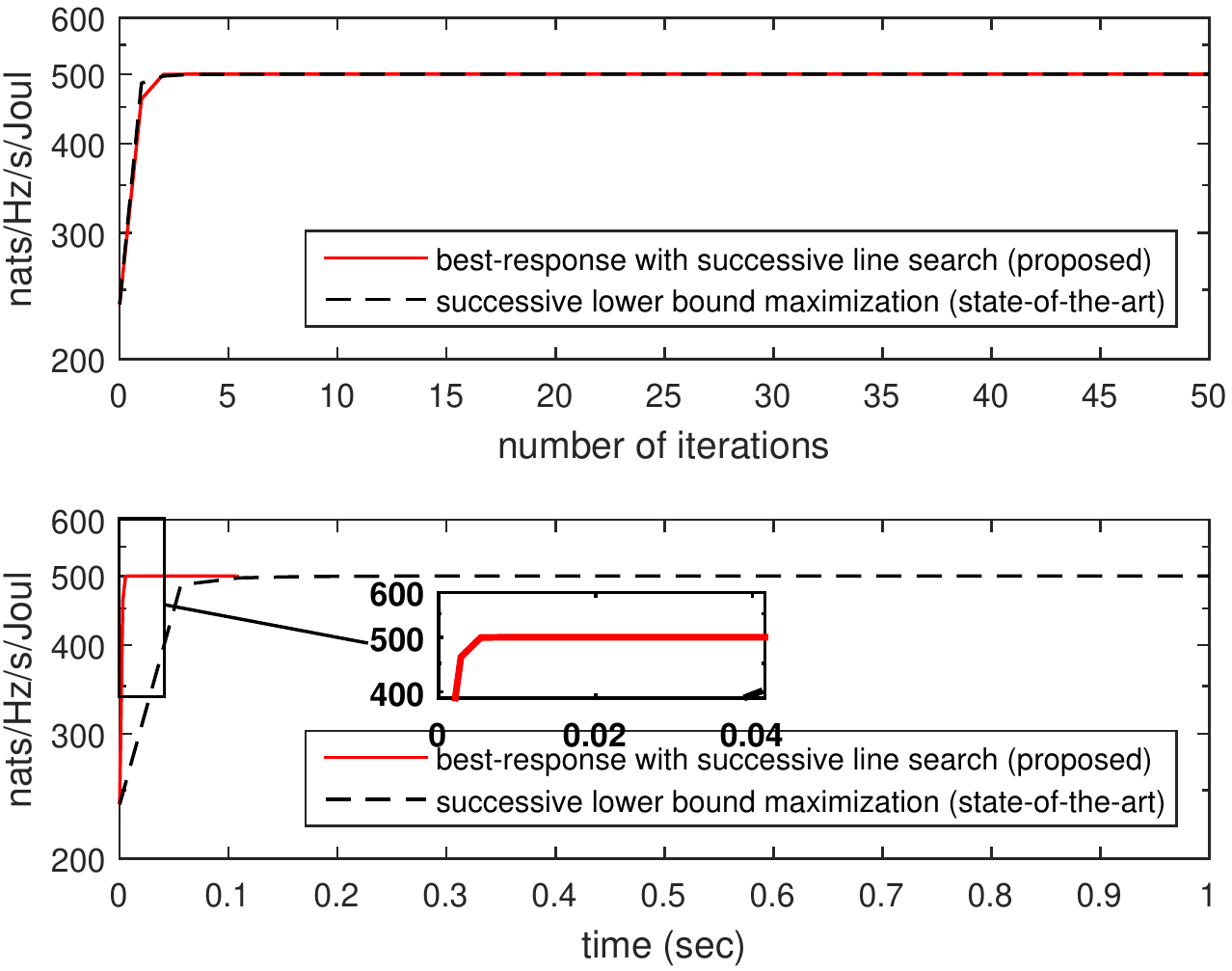}}\hspace{0.5cm}\subfigure[Number of users: $K=50$]{\includegraphics[scale=0.67]{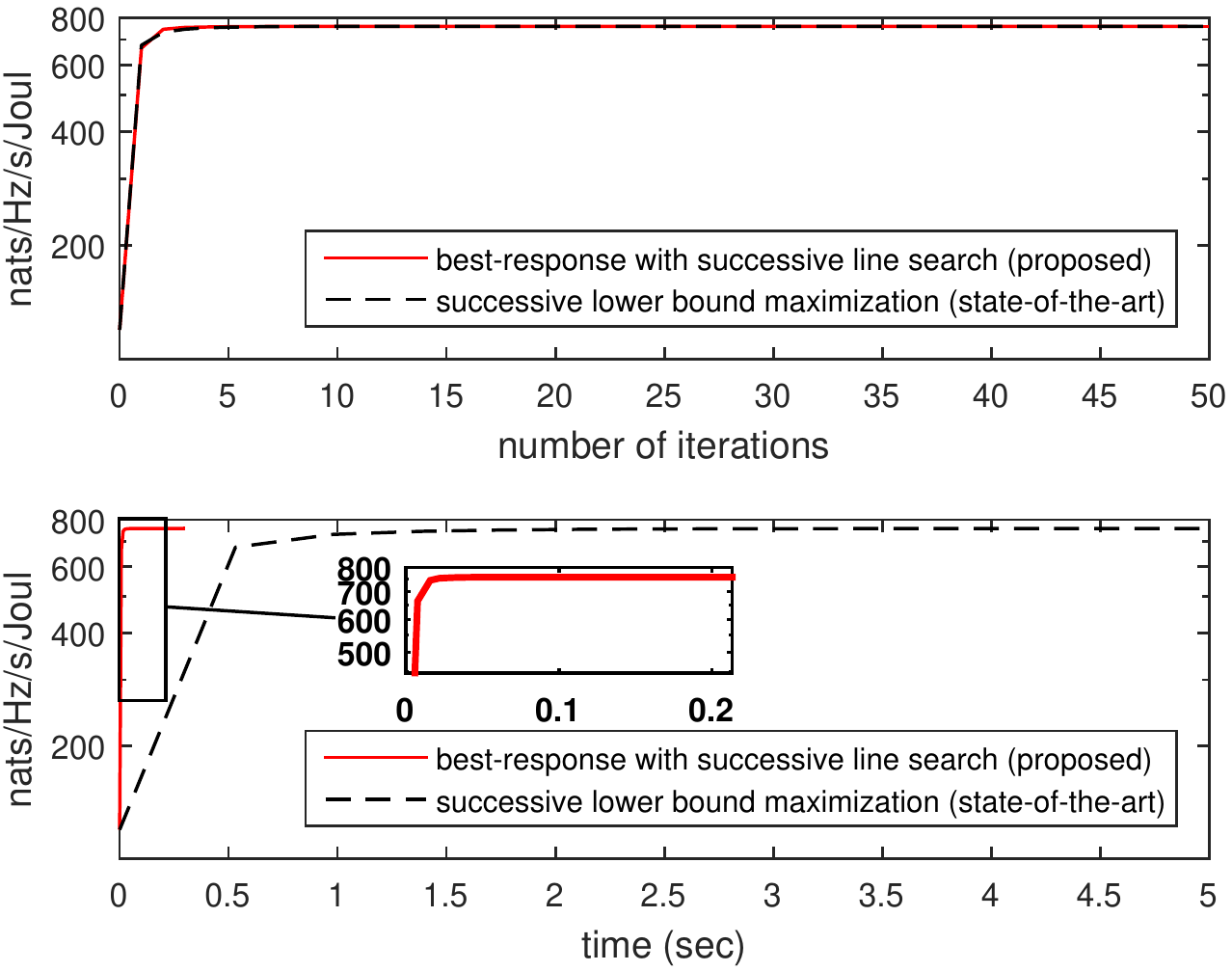}}\protect\caption{\label{fig:EE-convergence}Energy Efficiency Maximization: achieved
energy efficiency versus the number of iterations}
\end{figure*}

\textbf{Simulations.} The number of antennas at the BS in each cell
is $M=50$, and the channel from user $j$ to cell $k$ is $\mathbf{h}_{kj}\in\mathbb{C}^{M\times1}$.
We assume a similar setup as \cite{Zappone2015}: $w_{kk}=\left|\mathbf{h}_{kk}^{H}\mathbf{h}_{kk}\right|^{2}$,
$w_{kj}=\left|\mathbf{h}_{kk}^{H}\mathbf{h}_{kj}\right|^{2}+\epsilon\mathbf{h}_{kk}^{H}\mathbf{D}_{j}\mathbf{h}_{kk}$
for $j\neq k$ and $\phi_{k}=\epsilon\mathbf{h}_{kk}^{H}\mathbf{D}_{k}\mathbf{h}_{kk}$,
where $\epsilon=0.01$ is the error magnitude of hardware impairments
at the BS and $\mathbf{D}_{j}=\textrm{diag}(\{|h_{jj}(m)|^{2}\}_{m=1}^{M})$.
The noise covariance $\sigma_{k}^{2}=1$, and the hardware dissipated
power $p_{c}$ is 10dBm, while $\underline{p}_{k}$ is -10dBm and
$\overline{p}_{k}$ is 10dBm for all users. The benchmark algorithm
is \cite[Algorithm 1]{Zappone2015}, which successively maximizes
the following lower bound function of the objective function in (\ref{eq:EE-problem-formulation}),
which is tight at $\mathbf{p}=\mathbf{p}^{t}$:\addtocounter{equation}{1}
\begin{align}
\underset{\mathbf{q}}{\textrm{maximize}\;} & \frac{\sum_{k=1}^{K}b_{k}^{t}+a_{k}^{t}\log w_{kk}}{P_{c}+\sum_{k=1}^{K}e^{q_{k}}}+\nonumber \\
 & \frac{\sum_{k=1}^{K}a_{k}^{t}(q_{k}-\log(\sigma_{k}^{2}+\phi_{k}e^{q_{k}}+\sum_{j\neq k}w_{kj}e^{q_{j}}))}{P_{c}+\sum_{k=1}^{K}e^{q_{k}}}\nonumber \\
\textrm{subject to}\; & \log(\underline{p}_{k})\leq q_{k}\leq\log(\overline{p}_{k}),\, k=1,\ldots,K,\label{eq:EE-Jorswieck-ApproximateProblem}
\end{align}
where
\begin{align*}
a_{k}^{t} & \triangleq\frac{\textrm{sinr}_{k}(\mathbf{p}^{t})}{1+\textrm{sinr}_{k}(\mathbf{p}^{t})},\\
b_{k}^{t} & \triangleq\log(1+\textrm{sinr}_{k}(\mathbf{p}^{t}))-\frac{\textrm{sinr}_{k}(\mathbf{p}^{t})}{1+\textrm{sinr}_{k}(\mathbf{p}^{t})}\textrm{log}(\textrm{sinr}_{k}(\mathbf{p}^{t})),
\end{align*}
and
\[
\textrm{sinr}_{k}(\mathbf{p})\triangleq\frac{w_{kk}p^{t}}{\sigma_{k}^{2}+\phi_{k}p_{k}+\sum_{j\neq k}w_{kj}p_{j}}.
\]
Denote the optimal variable of (\ref{eq:EE-Jorswieck-ApproximateProblem})
as $\mathbf{q}^{t}$ (which can be found by the Dinkelbach's algorithm);
then the variable $\mathbf{p}$ is updated as $p_{k}^{t+1}=e^{q_{k}^{t}}$
for all $k=1,\ldots,K$. We thus coin \cite[Algorithm 1]{Zappone2015}
as the successive lower bound maximization (SLBM) method.

In Figure \ref{fig:EE-convergence}, we compare the convergence behavior
of the proposed method and the SLBM method in terms of both the number
of iterations (the upper subplots) and the CPU time (the lower subplots),
for two different number of users: $K=10$ in Figure \ref{fig:EE-convergence}
(a) and $K=50$ in Figure \ref{fig:EE-convergence} (b). It is obvious
that the convergence speed of the proposed algorithm in terms of the
number of iterations is comparable to that of the SLBM method. However,
we remark that the approximate problem (\ref{eq:EE-Approximate-problem-overall})
of the proposed algorithm is superior to that of the SLBM method in
the following aspects:

Firstly, the approximate problem of the proposed algorithm consists
of independent subproblems that can be solved in parallel, cf. (\ref{eq:EE-Approximate-problem-decomposition}),
while each subproblem has a closed-form solution, cf. (\ref{eq:EE-closed-form})-(\ref{eq:EE-closed-form-3}).
However, the optimization variable in the approximate problem of the
SLBM method (\ref{eq:EE-Jorswieck-ApproximateProblem}) is a vector
$\mathbf{q}\in\mathbb{R}^{K\times1}$ and the approximate problem
can only be solved by a general purpose solver.

In the simulations, we use the Matlab optimization toolbox to solve
(\ref{eq:EE-Jorswieck-ApproximateProblem}) and the iterative update
specified in (\ref{eq:EE-closed-form})-(\ref{eq:EE-closed-form-3})
to solve (\ref{eq:EE-Approximate-problem-overall}), where the stopping
criterion for (\ref{eq:EE-closed-form}) is $\left\Vert \boldsymbol{\lambda}^{t,\tau}\right\Vert _{\infty}\leq10^{-5}$.
The upper subplots in Figure \ref{fig:EE-convergence} show that the
numbers of iterations required for convergence are approximately the
same for the SLBM method when $K=10$ in Figure \ref{fig:EE-convergence}
(a) and when $K=50$ in Figure \ref{fig:EE-convergence} (b). However,
we see from the lower subplots in Figure \ref{fig:EE-convergence}
that the CPU time of each iteration of the SLBM method is dramatically
increased when $K$ is increased from 10 to 50. On the other hand,
the CPU time of the proposed algorithm is not notably changed because
the operations are parallelizable%
\footnote{By stacking the $p_{k}(\lambda_{k}^{t,\tau})$'s into the vector form
$\mathbf{p}(\boldsymbol{\lambda}^{t,\tau})=(p_{k}(\lambda_{k}^{t,\tau}))_{k=1}^{K}$
we can see that only element wise operations between vectors and matrix
vector multiplications are involved. The simulations on which Figure
\ref{fig:EE-convergence} are based are not performed in a real parallel
computing environment with $K$ processors, but only make use of the
efficient linear algebraic implementations available in Matlab which
already implicitly admits a certain level of parallelism.%
} and the required CPU time is thus not affected by the problem size.

Secondly, since a variable substitution $q_{k}=e^{p_{k}}$ is adopted
in the SLBM method (we refer to \cite{Zappone2015} for more details),
the lower bound constraint $\underline{p}_{k}=0$ (which corresponds
to $q_{k}=-\infty$) cannot be handled by the SLBM method numerically.
This limitation impairs the applicability of the SLBM method in many
practical scenarios.

\subsection{\label{sub:LASSO}LASSO}

In this subsection, we study the LASSO problem to illustrate the advantage
of the proposed line search method for nondifferentiable optimization
problems.

LASSO is an important and widely studied problem in sparse signal
recovery \cite{Tibshirani1996,Kim2007,Beck2009,Boyd2010}:
\begin{equation}
\begin{split}\underset{\mathbf{x}}{\textrm{minimize}}\quad & {\textstyle \frac{1}{2}}\left\Vert \mathbf{Ax-b}\right\Vert _{2}^{2}+\mu\left\Vert \mathbf{x}\right\Vert _{1},\end{split}
\label{eq:lasso}
\end{equation}
where $\mathbf{A}\in\mathbb{R}^{N\times K}$ (with $N\ll K$), $\mathbf{b}\in\mathbb{R}^{K\times1}$
and $\mu>0$ are given parameters. Problem (\ref{eq:lasso}) is an
instance of the general problem structure defined in (\ref{eq:original_function_nonsmooth})
with the following decomposition:
\begin{equation}
f(\mathbf{x})\triangleq{\textstyle \frac{1}{2}}\left\Vert \mathbf{Ax-b}\right\Vert _{2}^{2},\quad\textrm{and}\quad g(\mathbf{x})\triangleq\mu\left\Vert \mathbf{x}\right\Vert _{1}.\label{eq:lasso-f=000026g}
\end{equation}

Problem (\ref{eq:lasso}) is convex, but its objective function is
nondifferentiable and it does not have a closed-form solution. To
apply Algorithm \ref{alg:Successive-approximation-method-nonsmooth},
the scalar decomposition $\mathbf{x}=(x_{k})_{k=1}^{K}$ is adopted.
Recalling (\ref{eq:nonsmooth-approximate-problem-1}) and (\ref{eq:jacobi-approximate-function}),
the approximate problem is
\begin{equation}
\mathbb{B}\mathbf{x}^{t}=\underset{\mathbf{x}}{\arg\min}\;\bigl\{{\textstyle \sum_{k=1}^{K}}f(x_{k},\mathbf{x}_{-k}^{t})+g(\mathbf{x})\bigr\}.\label{eq:lasso-approximate}
\end{equation}
Note that $g(\mathbf{x})$ can be decomposed among different components
of $\mathbf{x}$, i.e., $g(\mathbf{x})=\sum_{k=1}^{K}g(x_{k})$, so
the vector problem (\ref{eq:lasso-approximate}) reduces to $K$ independent
scalar subproblems that can be solved in parallel:
\begin{align*}
\mathbb{B}_{k}\mathbf{x}^{t} & =\underset{x_{k}}{\arg\min}\;\bigl\{ f(x_{k},\mathbf{x}_{-k}^{t})+g(x_{k})\bigr\}\\
 & =d_{k}(\mathbf{A}^{T}\mathbf{A})^{-1}\mathcal{S}_{\mu}(r_{k}(\mathbf{x}^{t})),\, k=1,\ldots,K,
\end{align*}
where $d_{k}(\mathbf{A}^{T}\mathbf{A})$ is the $k$-th diagonal element
of $\mathbf{A}^{T}\mathbf{A}$, $\mathcal{S}_{\mathbf{a}}(\mathbf{b})\triangleq\left[\mathbf{b-a}\right]^{+}-\left[\mathbf{-b-a}\right]^{+}$
is the so-called soft-thresholding operator \cite{Beck2009} and
\begin{equation}
\mathbf{r}(\mathbf{x})\triangleq\mathbf{d}(\mathbf{A}^{T}\mathbf{A})\circ\mathbf{x}-\mathbf{A}^{T}(\mathbf{Ax-b}),\label{eq:r(x)}
\end{equation}
or more compactly:
\begin{equation}
\mathbb{B}\mathbf{x}^{t}=(\mathbb{B}_{k}\mathbf{x}^{t})_{k=1}^{K}=\mathbf{d}(\mathbf{A}^{T}\mathbf{A})^{-1}\circ\mathcal{S}_{\mu\mathbf{1}}(\mathbf{r}(\mathbf{x}^{t})).\label{eq:lasso-approximate-problem}
\end{equation}
Thus the update direction exhibits a closed-form expression. The stepsize
based on the proposed exact line search (\ref{eq:nonsmooth-stepsize-1})
is
\begin{align}
\gamma^{t} & =\underset{0\leq\gamma\leq1}{\arg\min}\;\left\{ \negthickspace\negthickspace\begin{array}{l}
f(\mathbf{x}^{t}+\gamma(\mathbb{B}\mathbf{x}^{t}-\mathbf{x}^{t}))+\gamma\bigl(g(\mathbb{B}\mathbf{x}^{t})-g(\mathbf{x}^{t})\bigr)\end{array}\negthickspace\negthickspace\right\} \nonumber \\
 & =\underset{0\leq\gamma\leq1}{\arg\min}\;\left\{ \negthickspace\negthickspace\begin{array}{l}
\frac{1}{2}\left\Vert \mathbf{A}(\mathbf{x}^{t}+\gamma(\mathbb{B}\mathbf{x}^{t}-\mathbf{x}^{t}))-\mathbf{b}\right\Vert _{2}^{2}\smallskip\\
\quad+\,\gamma\,\mu\bigl(\left\Vert \mathbb{B}\mathbf{x}^{t}\right\Vert _{1}-\left\Vert \mathbf{x}^{t}\right\Vert _{1}\bigr)
\end{array}\negthickspace\negthickspace\right\} \nonumber \\
 & =\left[-\frac{(\mathbf{A}\mathbf{x}^{t}-\mathbf{b})^{T}\mathbf{A}(\mathbb{B}\mathbf{x}^{t}-\mathbf{x}^{t})+\mu(\left\Vert \mathbb{B}\mathbf{x}^{t}\right\Vert _{1}-\left\Vert \mathbf{x}^{t}\right\Vert _{1})}{(\mathbf{A}(\mathbb{B}\mathbf{x}^{t}-\mathbf{x}^{t}))^{T}(\mathbf{A}(\mathbb{B}\mathbf{x}^{t}-\mathbf{x}^{t}))}\right]_{0}^{1}.\label{eq:lasso-approximate-stepsize}
\end{align}
The exact line search consists in solving a convex quadratic optimization
problem with a scalar variable and a bound constraint, so the problem
exhibits a closed-form solution (\ref{eq:lasso-approximate-stepsize}).
Therefore, both the update direction and stepsize can be calculated
in closed-form. We name the proposed update (\ref{eq:lasso-approximate-problem})-(\ref{eq:lasso-approximate-stepsize})
as Soft-Thresholding with Exact Line search Algorithm ($\texttt{STELA}$).

The proposed update (\ref{eq:lasso-approximate-problem})-(\ref{eq:lasso-approximate-stepsize})
has several desirable features that make it appealing in practice.
Firstly, in each iteration, all elements are updated in parallel based
on the nonlinear best-response (\ref{eq:lasso-approximate-problem}).
This is in the same spirit as \cite{Scutari_BigData,Elad2006} and
the convergence speed is generally faster than BCD \cite{Friedman2007}
or the gradient-based update \cite{Figueiredo2007}. Secondly, the
proposed exact line search (\ref{eq:lasso-approximate-stepsize})
not only yields notable progress in each iteration but also enjoys
an easy implementation given the closed-form expression. The convergence
speed is thus further enhanced as compared to the procedures proposed
in \cite{Tseng2009,Elad2006,Scutari_BigData} where either decreasing
stepsizes are used \cite{Scutari_BigData} or the line search is over
the original nondifferentiable objective function in (\ref{eq:lasso})
\cite{Elad2006,Tseng2009}:
\[
\min_{0\leq\gamma\leq1}\left\{ \begin{array}{l}
\frac{1}{2}\left\Vert \mathbf{A}(\mathbf{x}^{t}+\gamma(\mathbb{B}\mathbf{x}^{t}-\mathbf{x}^{t}))-\mathbf{b}\right\Vert _{2}^{2}\smallskip\\
+\mu\left\Vert \mathbf{x}^{t}+\gamma(\mathbb{B}\mathbf{x}^{t}-\mathbf{x}^{t})\right\Vert _{1}
\end{array}\right\} .
\]

\textbf{Computational complexity}. The computational overhead associated
with the proposed exact line search (\ref{eq:lasso-approximate-stepsize})
can significantly be reduced if (\ref{eq:lasso-approximate-stepsize})
is carefully implemented as outlined in the following. The most complex
operation in (\ref{eq:lasso-approximate-stepsize}) is the matrix-vector
multiplication, namely, $\mathbf{A}\mathbf{x}^{t}-\mathbf{b}$ in
the numerator and $\mathbf{A}(\mathbb{B}\mathbf{x}^{t}-\mathbf{x}^{t})$
in the denominator. On the one hand, the term $\mathbf{A}\mathbf{x}^{t}-\mathbf{b}$
is already available from $\mathbf{r}(\mathbf{x}^{t})$, which is
computed in order to determine the best-response in (\ref{eq:lasso-approximate-problem}).
On the other hand, the matrix-vector multiplication $\mathbf{A}(\mathbb{B}\mathbf{x}^{t}-\mathbf{x}^{t})$
is also required for the computation of $\mathbf{A}\mathbf{x}^{t+1}-\mathbf{b}$
as it can alternatively be computed as:
\begin{align}
\mathbf{A}\mathbf{x}^{t+1}-\mathbf{b} & =\mathbf{A}(\mathbf{x}^{t}+\gamma^{t}(\mathbb{B}\mathbf{x}^{t}-\mathbf{x}^{t}))-\mathbf{b}\nonumber \\
 & =(\mathbf{A}\mathbf{x}^{t}-\mathbf{b})+\gamma^{t}\mathbf{A}(\mathbb{B}\mathbf{x}^{t}-\mathbf{x}^{t}),\label{eq:s_(t+1)}
\end{align}
where then only an additional vector addition is involved. As a result,
the stepsize (\ref{eq:lasso-approximate-stepsize}) does not incur
any additional matrix-vector multiplications, but only affordable
vector-vector multiplications.

\begin{figure}[t]
\center\includegraphics[bb=12bp 367bp 555bp 825bp,clip,scale=0.45]{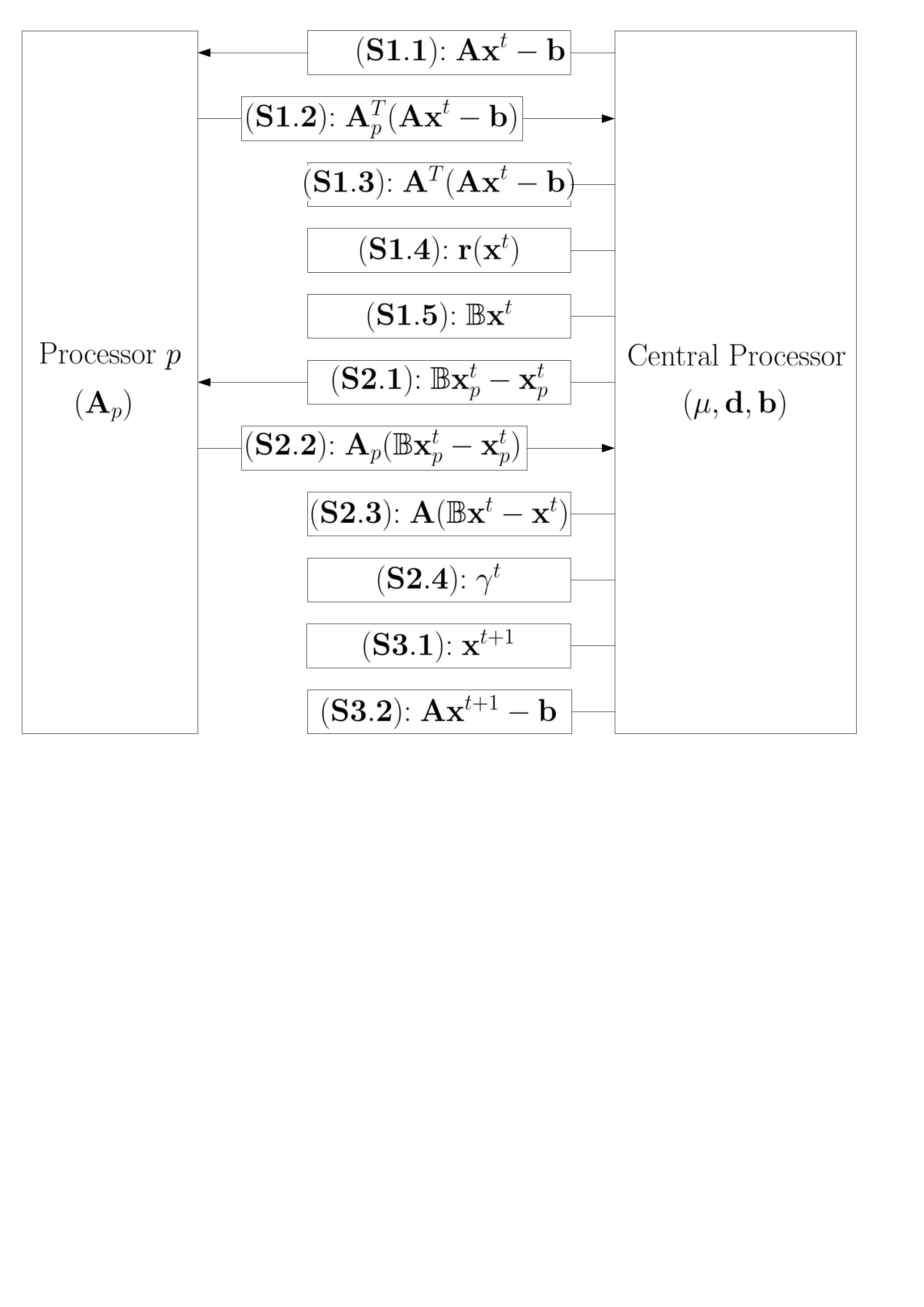}\protect\caption{\label{fig:lasso-Signaling-exchange}Operation flow and signaling
exchange between local processor $p$ and the central processor. A
solid line indicates the computation that is locally performed by
the central/local processor, and a solid line with an arrow indicates
signaling exchange between the central and local processor and the
direction of the signaling exchange.}

\vspace{-1em}
\end{figure}

\textbf{Signaling exchange}. When $\mathbf{A}$ cannot be stored and
processed by a centralized processing unit, a parallel architecture
can be employed. Assume there are $P+1$ ($P\geq2$) processors. We
label the first $P$ processors as local processors and the last one
as the central processor, and partition $\mathbf{A}$ as
\[
\mathbf{A}=[\mathbf{A}_{1},\;\mathbf{A}_{2},\;\ldots\;,\mathbf{A}_{P}],
\]
where $\mathbf{A}_{p}\in\mathbb{R}^{N\times K_{p}}$ and $\sum_{p=1}^{P}K_{p}=K$.
Matrix $\mathbf{A}_{p}$ is stored and processed in the local processor
$p$, and the following computations are decomposed among the local
processors:\begin{subequations}\label{eq:decompose}
\begin{align}
\mathbf{A}\mathbf{x} & ={\textstyle \sum_{p=1}^{P}}\mathbf{A}_{p}\mathbf{x}_{p},\label{eq:decompose-1}\\
\mathbf{A}^{T}(\mathbf{Ax-b}) & =\left(\mathbf{A}_{p}^{T}(\mathbf{Ax-b})\right)_{p=1}^{P},\label{eq:decompose-2}\\
\mathbf{d}(\mathbf{A}^{T}\mathbf{A}) & =(\mathbf{d}(\mathbf{A}_{p}^{T}\mathbf{A}_{p}))_{p=1}^{P}.\label{eq:decompose-3}
\end{align}
\end{subequations}where $\mathbf{x}_{p}\in\mathbb{R}^{K_{p}}$. The
central processor computes the best-response $\mathbb{B}\mathbf{x}^{t}$
in (\ref{eq:lasso-approximate-problem}) and the stepsize $\gamma^{t}$
in (\ref{eq:lasso-approximate-stepsize}). The decomposition in (\ref{eq:decompose})
enables us to analyze the signaling exchange between local processor
$p$ and the central processor involved in (\ref{eq:lasso-approximate-problem})
and (\ref{eq:lasso-approximate-stepsize})%
\footnote{Updates (\ref{eq:lasso-approximate-problem}) and (\ref{eq:lasso-approximate-stepsize})
can also be implemented by a parallel architecture without a central
processor. In this case, the signaling is exchanged mutually between
every two of the local processors, but the analysis is similar and
the conclusion to be drawn remains same: the proposed exact line search
(\ref{eq:lasso-approximate-stepsize}) does not incur additional signaling
compared with predetermined stepsizes.%
}.

The signaling exchange is summarized in Figure \ref{fig:lasso-Signaling-exchange}.
Firstly, the central processor sends $\mathbf{A}\mathbf{x}^{t}-\mathbf{b}$
to the local processors (\textbf{S1.1})%
\footnote{$\mathbf{x}^{0}$ is set to $\mathbf{x}^{0}=\mathbf{0}$, so $\mathbf{A}\mathbf{x}^{0}-\mathbf{b}=-\mathbf{b}$.%
}, and each local processor $p$ for $p=1,\ldots,P$ first computes
$\mathbf{A}_{p}^{T}(\mathbf{A}\mathbf{x}^{t}-\mathbf{b})$ and then
sends it back to the central processor (\textbf{S1.2}), which forms
$\mathbf{A}^{T}(\mathbf{A}\mathbf{x}^{t}-\mathbf{b})$ (\textbf{S1.3})
as in (\ref{eq:decompose-2}) and calculates $\mathbf{r}(\mathbf{x}^{t})$
as in (\ref{eq:r(x)}) (\textbf{S1.4}) and then $\mathbb{B}\mathbf{x}^{t}$
as in (\ref{eq:lasso-approximate-problem}) (\textbf{S1.5}). Then
the central processor sends $\mathbb{B}\mathbf{x}_{p}^{t}-\mathbf{x}_{p}^{t}$
to the local processor $p$ for $p=1,\ldots,P$ (\textbf{S2.1}), and
each local processor first computes $\mathbf{A}_{p}(\mathbb{B}\mathbf{x}_{p}^{t}-\mathbf{x}_{p}^{t})$
and then sends it back to the central processor (\textbf{S2.2}), which
forms $\mathbf{A}(\mathbb{B}\mathbf{x}^{t}-\mathbf{x}^{t})$ (\textbf{S2.3})
as in (\ref{eq:decompose-1}), calculates $\gamma^{t}$ as in (\ref{eq:lasso-approximate-stepsize})
(\textbf{S2.4}), and updates $\mathbf{x}^{t+1}$ (\textbf{S3.1}) and
$\mathbf{A}\mathbf{x}^{t+1}-\mathbf{b}$ (\textbf{S3.2}) according
to (\ref{eq:s_(t+1)}). From Figure \ref{fig:lasso-Signaling-exchange}
we observe that the exact line search \eqref{eq:lasso-approximate-stepsize}
does \emph{not} incur any additional signaling compared with that
of predetermined stepsizes (e.g., constant and decreasing stepsize),
because the signaling exchange in \textbf{S2.1-S2.2 }has also to be
carried out in the computation of $\mathbf{A}\mathbf{x}^{t+1}-\mathbf{b}$
in \textbf{S3.2}, cf. (\ref{eq:s_(t+1)}).

We finally remark that the proposed successive line search can also
be applied and it exhibits a closed-form expression as well. However,
since the exact line search yields faster convergence, we omit the
details at this point.

\begin{figure}[t]
\center\includegraphics[scale=0.63]{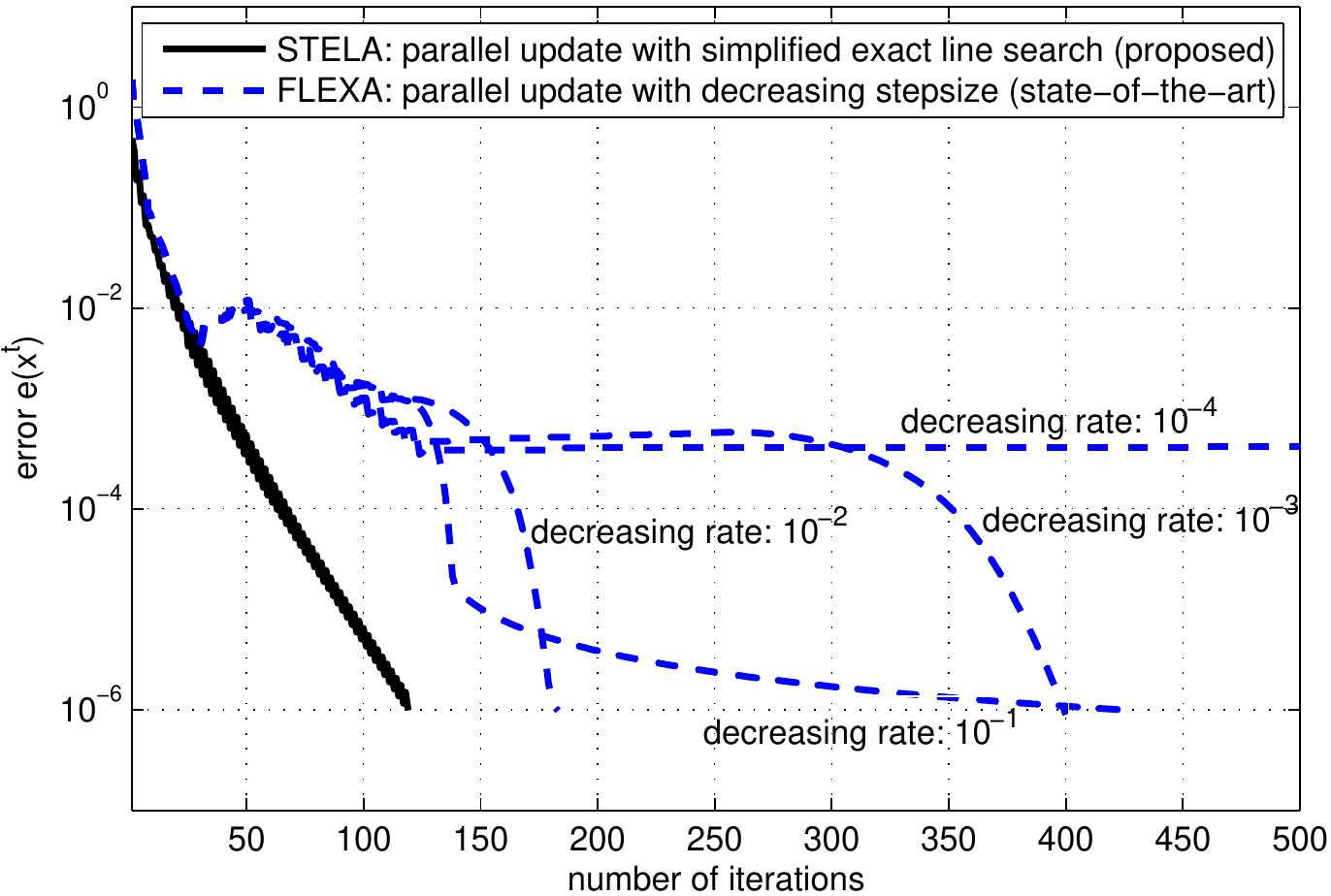}

\protect\caption{\label{fig:LASSO_decreasing_fails}Convergence of $\texttt{STELA}$
(proposed) and $\texttt{FLEXA}$ (state-of-the-art) for LASSO: error
versus the number of iterations.}
\end{figure}

\textbf{Simulations}.\textbf{ }We first compare in Figure \ref{fig:LASSO_decreasing_fails}
the proposed algorithm $\texttt{STELA}$ with $\texttt{FLEXA}$ \cite{Scutari_BigData}
in terms of the error criterion $e(\mathbf{x}^{t})$ defined as:
\begin{equation}
e(\mathbf{x}^{t})\triangleq\bigl\Vert\nabla f(\mathbf{x}^{t})-\left[\nabla f(\mathbf{x}^{t})-\mathbf{x}^{t}\right]_{-\mu\mathbf{1}}^{\mu\mathbf{1}}\bigr\Vert_{2}.\label{eq:lasso-error}
\end{equation}
Note that $\mathbf{x}^{\star}$ is a solution of (\ref{eq:lasso})
if and only if $e(\mathbf{x}^{\star})=0$ \cite{Byrd2013}. $\texttt{FLEXA}$
is implemented as outlined in \cite{Scutari_BigData}; however, the
selective update scheme \cite{Scutari_BigData} is not implemented
in $\texttt{FLEXA}$ because it is also applicable for $\texttt{STELA}$
and it cannot eliminate the slow convergence and sensitivity of the
decreasing stepsize. We also remark that the stepsize rule for $\texttt{FLEXA}$
is $\gamma^{t+1}=\gamma^{t}(1-\min(1,10^{-4}/e(\mathbf{x}^{t}))d\gamma^{t})$
\cite{Scutari_BigData}, where $d$ is the decreasing rate and $\gamma^{0}=0.9$.
The code and the data generating the figure can be downloaded online
\cite{Yang_code}.

Note that the error $e(\mathbf{x}^{t})$ plotted in Figure \ref{fig:LASSO_decreasing_fails}
does not necessarily decrease monotonically while the objective function
$f(\mathbf{x}^{t})+g(\mathbf{x}^{t})$ always does. This is because
$\texttt{STELA}$ and $\texttt{FLEXA}$ are descent direction methods.
For $\texttt{FLEXA}$, when the decreasing rate is low ($d=10^{-4}$),
no improvement is observed after 100 iterations. As a matter of fact,
the stepsize in those iterations is so large that the function value
is actually dramatically increased, and thus the associated iterations
are discarded in Figure \ref{fig:LASSO_decreasing_fails}. A similar
behavior is also observed for $d=10^{-3}$, until the stepsize becomes
sufficiently small. When the stepsize is quickly decreasing ($d=10^{-1}$),
although improvement is made in all iterations, the asymptotic convergence
speed is slow because the stepsize is too small to make notable improvement.
For this example, the choice $d=10^{-2}$ performs well, but the value
of a good decreasing rate depends on the parameter setup (e.g., $\mathbf{A}$,
$\mathbf{b}$ and $\mu$) and no general rule performs equally well
for all choices of parameters. By comparison, the proposed algorithm
$\texttt{STELA}$ is fast to converge and exhibits stable performance
without requiring any parameter tuning.

\begin{figure}[t]
\center\includegraphics[scale=0.75]{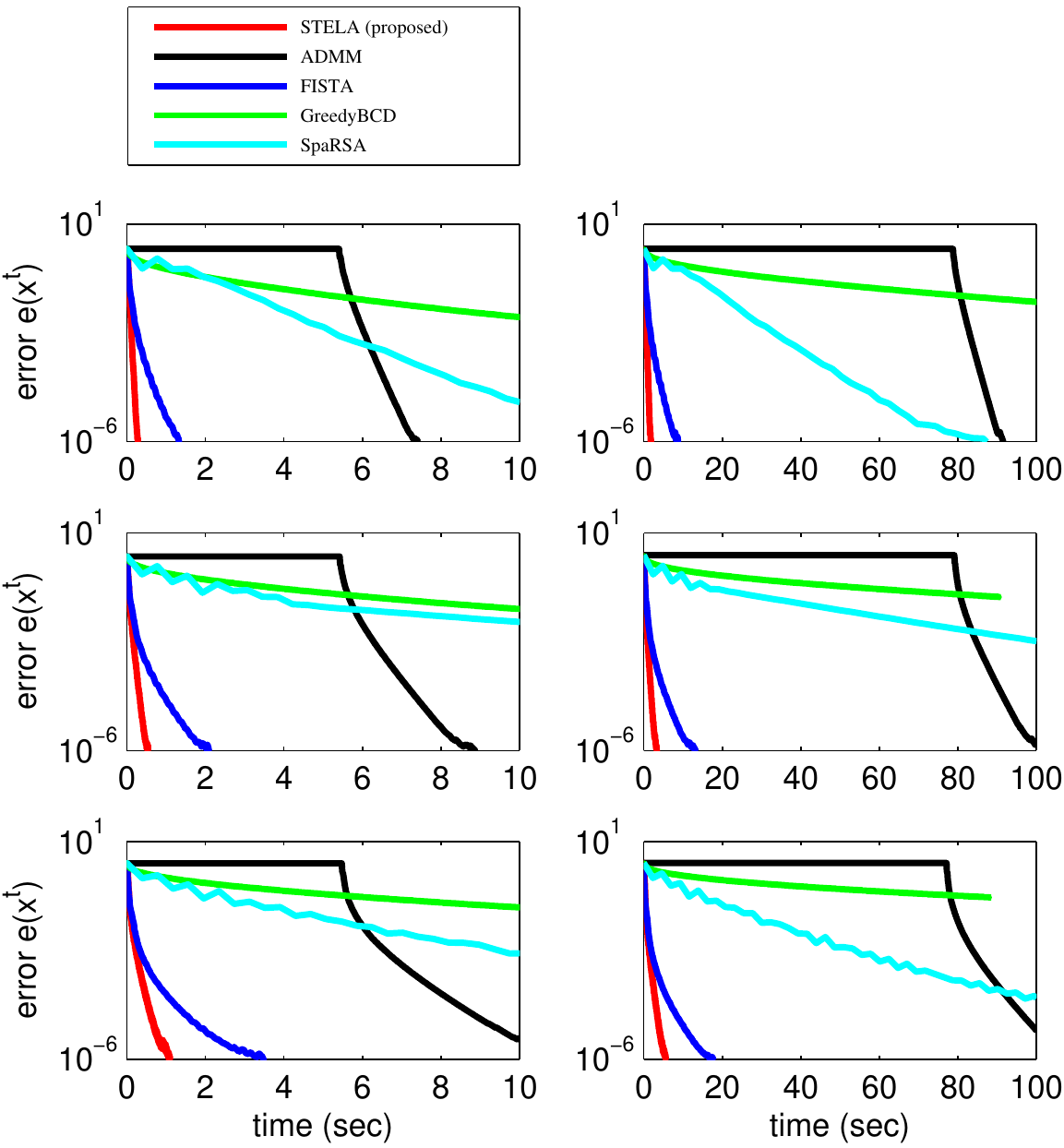}\protect\caption{\label{fig:LASSO_other_algorithms}Time versus error of different
algorithms for LASSO. In the left and right column, the dimension
of $\mathbf{A}$ is $2000\times4000$ and $5000\times10000$, respectively.
In the higher, middle and lower column, the density of $\mathbf{x}_{\textrm{true}}$
is 0.1, 0.2 and 0.4.}
\end{figure}

We also compare in Figure \ref{fig:LASSO_other_algorithms} the proposed
algorithm $\texttt{STELA}$ with other competitive algorithms in literature:
$\texttt{FISTA}$ \cite{Beck2009}, $\texttt{ADMM}$ \cite{Boyd2010},
$\texttt{GreedyBCD}$ \cite{Peng2013} and $\texttt{SpaRSA}$ \cite{Wright2009}.
We simulated $\texttt{GreedyBCD}$ of \cite{Peng2013} because it
exhibits guaranteed convergence. The dimension of $\mathbf{A}$ is
$2000\times4000$ (the left column of Figure \ref{fig:LASSO_other_algorithms})
and $5000\times10000$ (the right column). It is generated by the
Matlab command $\texttt{randn}$ with each row being normalized to
unity. The density (the proportion of nonzero elements) of the sparse
vector $\mathbf{x}_{\textrm{true}}$ is 0.1 (the upper row of Figure
\ref{fig:LASSO_other_algorithms}), 0.2 (the middle row) and 0.4 (the
lower row). The vector $\mathbf{b}$ is generated as $\mathbf{b}=\mathbf{A}\mathbf{x}_{\textrm{true}}+\mathbf{e}$
where $\mathbf{e}$ is drawn from an i.i.d. Gaussian distribution
with variance $10^{-4}$. The regularization gain $\mu$ is set to
$\mu=0.1\left\Vert \mathbf{A}^{T}\mathbf{b}\right\Vert _{\infty}$,
which allows $\mathbf{x}_{\textrm{true}}$ to be recovered to a high
accuracy \cite{Wright2009}.

The simulations are carried out under Matlab R2012a on a PC equipped
with an operating system of Windows 7 64-bit Home Premium Edition,
an Intel i5-3210 2.50GHz CPU, and a 8GB RAM. All of the Matlab codes
are available online \cite{Yang_code}. The comparison is made in
terms of CPU time that is required until either a given error bound
$e(\mathbf{x}^{t})\leq10^{-6}$ is reached or the maximum number of
iterations, namely, 2000, is reached. The running time consists of
both the initialization stage required for preprocessing (represented
by a flat curve) and the formal stage in which the iterations are
carried out. For example, in the proposed algorithm $\texttt{STELA}$,
$\mathbf{d}(\mathbf{A}^{T}\mathbf{A})$ is computed%
\footnote{The Matlab command is $\texttt{sum(A.\textasciicircum2,1)}$, so matrix-matrix
multiplication between $\mathbf{A}^{T}$ and $\mathbf{A}$ is not
required.%
} in the initialization stage since it is required in the iterative
variable update in the formal stage, cf. (\ref{eq:lasso-approximate-problem}).
The simulation results are averaged over 20 instances.

We observe from Figure \ref{fig:LASSO_other_algorithms} that the
proposed algorithm $\texttt{STELA}$ converges faster than all competing
algorithms. Some further observations are in order.

$\bullet$ The proposed algorithm $\texttt{STELA}$ is not sensitive
to the density of the true signal $\mathbf{x}_{\textrm{true}}$. When
the density is increased from 0.1 (left column) to 0.2 (middle column)
and then to 0.4 (right column), the CPU time increases negligibly.

$\bullet$ The proposed algorithm $\texttt{STELA}$ scales relatively
well with the problem dimension. When the dimension of $\mathbf{A}$
is increased from $2000\times4000$ (the left column) to $5000\times10000$
(the right column), the CPU time is only marginally increased.

$\bullet$ The initialization stage of $\texttt{ADMM}$ is time consuming
because of some expensive matrix operations as, e.g., $\mathbf{A}\mathbf{A}^{T}$,
$\left(\mathbf{I}+\frac{1}{c}\mathbf{A}\mathbf{A}^{T}\right)^{-1}$
and $\mathbf{A}^{T}\left(\mathbf{I}+\frac{1}{c}\mathbf{A}\mathbf{A}^{T}\right)^{-1}\mathbf{A}$
($c$ is a given positive constant). More details can be found in
\cite[Sec. 6.4]{Boyd2010}. Furthermore, the CPU time of the initialization
stage of $\texttt{ADMM}$ is increased dramatically when the dimension
of $\mathbf{A}$ is increased from $2000\times4000$ to $5000\times10000$.

$\bullet$ $\texttt{SpaRSA}$ performs better when the density of
$\mathbf{x}_{\textrm{true}}$ is smaller, e.g., 0.1, than in the case
when it is large, e.g., 0.2 and 0.4.

$\bullet$ The asymptotic convergence speed of $\texttt{GreedyBCD}$
is slow, because only one variable is updated in each iteration.

To further evaluate the performance of the proposed algorithm $\texttt{STELA}$,
we test it on the benchmarking platform developed by the Optimization
Group from the Department of Mathematics at the Darmstadt University
of Technology%
\footnote{Project website: http://wwwopt.mathematik.tu-darmstadt.de/spear/%
} and compare it with different algorithms in various setups (data
set, problem dimension, etc.) for the basis pursuit (BP) problem \cite{Lorenz2015}:
\[
\begin{split}\textrm{minimize}\quad & \left\Vert \mathbf{x}\right\Vert _{1}\\
\textrm{subject to}\quad & \mathbf{Ax=b}.
\end{split}
\]
To adapt $\texttt{STELA}$for the BP problem, we use the augmented
Lagrangian approach \cite{Rockafellar1976a,bertsekas1999nonlinear}:
\begin{eqnarray*}
\mathbf{x}^{t+1} & = & \left\Vert \mathbf{x}\right\Vert _{1}+(\boldsymbol{\lambda}^{t})^{T}(\mathbf{Ax-b})+\frac{c^{t}}{2}\left\Vert \mathbf{Ax-b}\right\Vert _{2}^{2},\\
\boldsymbol{\lambda}^{t+1} & = & \boldsymbol{\lambda^{t}}+c^{t}(\mathbf{A}\mathbf{x}^{t+1}-\mathbf{b}),
\end{eqnarray*}
where $c^{t+1}=\min(2c^{t},10^{2})$ ($c^{0}=10/\left\Vert \mathbf{A}^{T}\mathbf{b}\right\Vert _{\infty}$),
$\mathbf{x}^{t+1}$ is computed by $\texttt{STELA}$ and this process
is repeated until $\boldsymbol{\lambda}^{t}$ converges. The numerical
results summarized in \cite{JanTillmann} show that, although $\texttt{STELA}$
must be called multiple times before the Lagrange multiplier $\boldsymbol{\lambda}$
converges, the proposed algorithm for BP based on $\texttt{STELA}$
is very competitive in terms of running time and robust in the sense
that it solved all problem instances in the test platform database.

\section{\label{sec:Concluding-remarks}Concluding Remarks}

In this paper, we have proposed a novel iterative algorithm based
on convex approximation. The most critical requirement on the approximate
function is that it is pseudo-convex. On the one hand, the relaxation
of the assumptions on the approximate functions can make the approximate
problems much easier to solve. We show by a counter-example that the
assumption on pseudo-convexity is tight in the sense that when it
is violated, the algorithm may not converge. On the another hand,
the stepsize based on the exact/successive line search yields notable
progress in each iteration. Additional structures can be exploited
to assist with the selection of the stepsize, so that the algorithm
can be further accelerated. The advantages and benefits of the proposed
algorithm have been demonstrated using prominent applications in communication
networks and signal processing, and they are also numerically consolidated.
The proposed algorithm can readily be applied to solve other problems
as well, such as portfolio optimization \cite{Yang2013b}.

\appendices{}

\section{\label{sec:Proof-of-Proposition-descent}Proof of Proposition \ref{prop:descent-property}}
\begin{IEEEproof}
i) Firstly, suppose $\mathbf{y}$ is a stationary point of (\ref{eq:original_function});
it satisfies the first-order optimality condition:
\[
\nabla f(\mathbf{y})^{T}(\mathbf{x}-\mathbf{y})\geq0,\;\forall\,\mathbf{x}\in\mathcal{X}.
\]
Using Assumption (A3), we get
\[
\nabla\tilde{f}(\mathbf{y};\mathbf{y})^{T}(\mathbf{x}-\mathbf{y})\geq0,\;\forall\,\mathbf{x}\in\mathcal{X}.
\]
Since $\tilde{f}(\bullet;\mathbf{y})$ is pseudo-convex, the above
condition implies
\[
\tilde{f}(\mathbf{x};\mathbf{y})\geq\tilde{f}(\mathbf{y};\mathbf{y}),\;\forall\,\mathbf{x}\in\mathcal{X}.
\]
That is, $\tilde{f}(\mathbf{y};\mathbf{y})=\min_{\mathbf{x}\in\mathcal{X}}\tilde{f}(\mathbf{x};\mathbf{y})$
and $\mathbf{y}\in\mathcal{S}(\mathbf{y})$.

Secondly, suppose $\mathbf{y}\in\mathcal{S}(\mathbf{y})$. We readily
get
\begin{equation}
\nabla f(\mathbf{y})^{T}(\mathbf{x}-\mathbf{y})=\nabla\tilde{f}(\mathbf{y};\mathbf{y})^{T}(\mathbf{x}-\mathbf{y})\geq0,\;\forall\,\mathbf{x}\in\mathcal{X},\label{eq:minimum-principle-1}
\end{equation}
where the equality and inequality comes from Assumption (A3) and the
first-order optimality condition, respectively, so $\mathbf{y}$ is
a stationary point of (\ref{eq:original_function}).

ii) From the definition of $\mathbb{B}\mathbf{x}$, it is either\begin{subequations}
\begin{equation}
\tilde{f}(\mathbb{B}\mathbf{y};\mathbf{y})=\tilde{f}(\mathbf{y};\mathbf{y}),\label{eq:possibility-1}
\end{equation}
or
\begin{equation}
\tilde{f}(\mathbb{B}\mathbf{y};\mathbf{y})<\tilde{f}(\mathbf{y};\mathbf{y}),\label{eq:possibility-2}
\end{equation}
\end{subequations}If (\ref{eq:possibility-1}) is true, then $\mathbf{y}\in\mathcal{S}(\mathbf{y})$
and, as we have just shown, it is a stationary point of (\ref{eq:original_function}).
So only (\ref{eq:possibility-2}) can be true. We know from the pseudo-convexity
of $\tilde{f}(\mathbf{x};\mathbf{y})$ in $\mathbf{x}$ (cf. Assumption
(A1)) and (\ref{eq:possibility-2}) that $\mathbb{B}\mathbf{y}\neq\mathbf{y}$
and
\begin{equation}
\nabla\tilde{f}(\mathbf{y};\mathbf{y})^{T}(\mathbb{B}\mathbf{y}-\mathbf{y})=\nabla f(\mathbf{y})^{T}(\mathbb{B}\mathbf{y}-\mathbf{y})<0,\label{eq:descent-direction-1}
\end{equation}
where the equality comes from Assumption (A3).
\end{IEEEproof}

\section{\label{appendix:Proof-of-Theorem}Proof of Theorem \ref{thm:convergence}}
\begin{IEEEproof}
Since $\mathbb{B}\mathbf{x}^{t}$ is the optimal point of (\ref{eq:approximate-problem}),
it satisfies the first-order optimality condition:
\begin{equation}
\nabla\tilde{f}(\mathbb{B}\mathbf{x}^{t};\mathbf{x}^{t})^{T}(\mathbf{x}-\mathbb{B}\mathbf{x}^{t})\geq0,\;\forall\,\mathbf{x}\in\mathcal{X}.\label{eq:minimum-principle}
\end{equation}

If (\ref{eq:possibility-1}) is true, then $\mathbf{x}^{t}\in\mathcal{S}(\mathbf{x}^{t})$
and it is a stationary point of (\ref{eq:original_function}) according
to Proposition \ref{prop:descent-property} (i). Besides, it follows
from (\ref{eq:minimum-principle-1}) (with $\mathbf{x}=\mathbb{B}\mathbf{x}^{t}$
and $\mathbf{y}=\mathbf{x}^{t}$) that $\nabla f(\mathbf{x}^{t})^{T}(\mathbb{B}\mathbf{x}^{t}-\mathbf{x}^{t})\geq0$.
Note that equality is actually achieved, i.e.,
\[
\nabla f(\mathbf{x}^{t})^{T}(\mathbb{B}\mathbf{x}^{t}-\mathbf{x}^{t})=0
\]
because otherwise $\mathbb{B}\mathbf{x}^{t}-\mathbf{x}^{t}$ would
be an ascent direction of $\tilde{f}(\mathbf{x};\mathbf{x}^{t})$
at $\mathbf{x}=\mathbf{x}^{t}$ and the definition of $\mathbb{B}\mathbf{x}^{t}$
would be contradicted. Then from the definition of the successive
line search, we can readily infer that
\begin{equation}
f(\mathbf{x}^{t+1})\leq f(\mathbf{x}^{t}).\label{eq:decreasing-1}
\end{equation}
It is easy to see (\ref{eq:decreasing-1}) holds for the exact line
search as well.

If (\ref{eq:possibility-2}) is true, $\mathbf{x}^{t}$ is not a stationary
point and $\mathbb{B}\mathbf{x}^{t}-\mathbf{x}^{t}$ is a strict descent
direction of $f(\mathbf{x})$ at $\mathbf{x}=\mathbf{x}^{t}$ according
to Proposition \ref{prop:descent-property} (ii): $f(\mathbf{x})$
is strictly decreased compared with $f(\mathbf{x}^{t})$ if $\mathbf{x}$
is updated at $\mathbf{x}^{t}$ along the direction $\mathbb{B}\mathbf{x}^{t}-\mathbf{x}^{t}$.
From the definition of the successive line search, there always exists
a $\gamma^{t}$ such that $0<\gamma^{t}\leq1$ and
\begin{equation}
f(\mathbf{x}^{t+1})=f(\mathbf{x}^{t}+\gamma^{t}(\mathbb{B}\mathbf{x}^{t}-\mathbf{x}^{t}))<f(\mathbf{x}^{t}).\label{eq:decreasing-2}
\end{equation}
This strict decreasing property also holds for the exact line search
because it is the stepsize that yields the largest decrease, which
is always larger than or equal to that of the successive line search.

We know from (\ref{eq:decreasing-1}) and (\ref{eq:decreasing-2})
that $\left\{ f(\mathbf{x}^{t})\right\} $ is a monotonically decreasing
sequence and it thus converges. Besides, for any two (possibly different)
convergent subsequences $\left\{ \mathbf{x}^{t}\right\} _{t\in\mathcal{T}_{1}}$
and $\left\{ \mathbf{x}^{t}\right\} _{t\in\mathcal{T}_{2}}$, the
following holds:
\[
\lim_{t\rightarrow\infty}f(\mathbf{x}^{t})=\lim_{\mathcal{T}_{1}\ni t\rightarrow\infty}f(\mathbf{x}^{t})=\lim_{\mathcal{T}_{2}\ni t\rightarrow\infty}f(\mathbf{x}^{t}).
\]
Since $f(\mathbf{x})$ is a continuous function, we infer from the
preceding equation that
\begin{equation}
f\left(\lim_{\mathcal{T}_{1}\ni t\rightarrow\infty}\mathbf{x}^{t}\right)=f\left(\lim_{\mathcal{T}_{2}\ni t\rightarrow\infty}\mathbf{x}^{t}\right).\label{eq:value-convergence}
\end{equation}

Now consider any convergent subsequence $\{\mathbf{x}^{t}\}_{t\in\mathcal{T}}$
with limit point $\mathbf{y}$, i.e., $\lim_{\mathcal{T}\ni t\rightarrow\infty}\mathbf{x}^{t}=\mathbf{y}$.
To show that $\mathbf{y}$ is a stationary point, we first assume
the contrary: $\mathbf{y}$ is not a stationary point. Since $\tilde{f}(\mathbf{x};\mathbf{x}^{t})$
is continuous in both $\mathbf{x}$ and $\mathbf{x}^{t}$ by Assumption
(A2) and $\left\{ \mathbf{B}\mathbf{x}^{t}\right\} _{t\in\mathcal{T}}$
is bounded by Assumption (A5), it follows from \cite[Th. 1]{Robinson1974}
that there exists a sequence $\left\{ \mathbb{B}\mathbf{x}^{t}\right\} _{t\in\mathcal{T}_{s}}$
with $\mathcal{T}_{s}\subseteq\mathcal{T}$ such that it converges
and $\lim_{\mathcal{T}_{s}\ni t\rightarrow\infty}\mathbb{B}\mathbf{x}^{t}\in\mathcal{S}(\mathbf{y})$.
Since both $f(\mathbf{x})$ and $\nabla f(\mathbf{x})$ are continuous,
applying \cite[Th. 1]{Robinson1974} again implies there is a $\mathcal{T}_{s'}$
such that $\mathcal{T}_{s'}\subseteq\mathcal{T}_{s}(\subseteq\mathcal{T})$
and $\left\{ \mathbf{x}^{t+1}\right\} _{t\in\mathcal{T}_{s'}}$ converges
to $\mathbf{y}'$ defined as:
\[
\mathbf{y}'\triangleq\mathbf{y}+\rho(\mathbb{B}\mathbf{y}-\mathbf{y}),
\]
where $\rho$ is the stepsize when either the exact or successive
line search is applied to $f(\mathbf{y})$ along the direction $\mathbb{B}\mathbf{y}-\mathbf{y}$.
Since $\mathbf{y}$ is not a stationary point, it follows from (\ref{eq:decreasing-2})
that $f(\mathbf{y}')<f(\mathbf{y})$, but this would contradict (\ref{eq:value-convergence}).
Therefore $\mathbf{y}$ is a stationary point, and the proof is completed.
\end{IEEEproof}

\section{\label{sec:Appendix-of-Theorem-Cartesian}Proof of Theorem \ref{thm:Cartesian}}
\begin{IEEEproof}
We first need to show that Proposition \ref{prop:descent-property}
still holds.

(i) We prove $\mathbf{y}$ is a stationary point of (\ref{eq:original-function-cartesian})
if and only if $\mathbf{y}_{k}\in\arg\min_{\mathbf{x}_{k}\in\mathcal{X}_{k}}f(\mathbf{x}_{k},\mathbf{y}_{-k})$
for all $k$.

Suppose $\mathbf{y}$ is a stationary point of (\ref{eq:original-function-cartesian}),
it satisfies the first-order optimality condition:
\[
\nabla f(\mathbf{y})^{T}(\mathbf{x}-\mathbf{y})={\textstyle \sum_{k=1}^{K}}\nabla_{k}f(\mathbf{y})^{T}(\mathbf{x}_{k}-\mathbf{y}_{k})\geq0,\forall\,\mathbf{x}\in\mathcal{X},
\]
and it is equivalent to
\[
\nabla_{k}f(\mathbf{y})^{T}(\mathbf{x}_{k}-\mathbf{y}_{k})\geq0,\forall\,\mathbf{x}_{k}\in\mathcal{X}_{k}.
\]
Since $f(\mathbf{x})$ is pseudo-convex in $\mathbf{x}_{k}$, the
above condition implies $f(\mathbf{y}_{k},\mathbf{y}_{-k})=\min_{\mathbf{x}_{k}\in\mathcal{X}_{k}}f(\mathbf{x}_{k},\mathbf{y}_{-k})$
for all $k=1,\ldots,K$.

Suppose $\mathbf{y}_{k}\in\arg\min_{\mathbf{x}_{k}\in\mathcal{X}_{k}}f(\mathbf{x}_{k},\mathbf{y}_{-k})$
for all $k=1,\ldots,K$. The first-order optimality conditions yields
\[
\nabla_{k}f(\mathbf{y})^{T}(\mathbf{x}_{k}-\mathbf{y}_{k})\geq0,\forall\,\mathbf{x}_{k}\in\mathcal{X}_{k}.
\]
Adding the above inequality for all $k=1,\ldots,K$ yields
\[
\nabla f(\mathbf{y})^{T}(\mathbf{x}-\mathbf{y})\geq0,\forall\,\mathbf{x}\in\mathcal{X}.
\]
Therefore, $\mathbf{y}$ is a stationary point of (\ref{eq:original-function-cartesian}).

(ii) We prove that if $\mathbf{y}$ is not a stationary point of (\ref{eq:original-function-cartesian}),
then $\nabla f(\mathbf{y})^{T}(\mathbb{B}\mathbf{y}-\mathbf{y})<0$.

It follows from the optimality of $\mathbb{B}_{k}\mathbf{x}$ that
\[
f(\mathbb{B}_{k}\mathbf{y},\mathbf{y}_{-k})\leq f(\mathbf{y}_{k},\mathbf{y}_{-k}),
\]
and
\begin{equation}
\nabla_{k}f(\mathbb{B}_{k}\mathbf{y},\mathbf{y}_{-k})^{T}(\mathbf{x}_{k}-\mathbb{B}_{k}\mathbf{y})\geq0,\forall\,\mathbf{x}_{k}\in\mathcal{X}_{k}.\label{eq:cartesian-proof-1}
\end{equation}

Firstly, there must exist an index $j$ such that
\begin{equation}
f(\mathbb{B}_{j}\mathbf{y},\mathbf{y}_{-j})<f(\mathbf{y}_{j},\mathbf{y}_{-j}),\label{eq:cartesian-proof-2}
\end{equation}
otherwise $\mathbf{y}$ would be a stationary point of (\ref{eq:original-function-cartesian}).
Since $f(\mathbf{x})$ is pseudo-convex in $\mathbf{x}_{k}$ for $k=1,\ldots,K$,
it follows from (\ref{eq:cartesian-proof-2}) that
\begin{equation}
\nabla_{j}f(\mathbf{y})^{T}(\mathbb{B}_{j}\mathbf{y}-\mathbf{y}_{j})<0.\label{eq:cartesian-proof-3}
\end{equation}

Secondly, for any index $k$ such that $f(\mathbb{B}_{k}\mathbf{y},\mathbf{y}_{-k})=f(\mathbf{y}_{k},\mathbf{y}_{-k})$,
$\mathbf{y}_{k}$ minimizes $f(\mathbf{x}_{k},\mathbf{y}_{-k})$ over
$\mathbf{x}_{k}\in\mathcal{X}_{k}$ and $\nabla_{k}f(\mathbf{y}_{k},\mathbf{y}_{-k})^{T}(\mathbf{x}_{k}-\mathbf{y}_{k})\geq0$
for any $\mathbf{x}_{k}\in\mathcal{X}$. Setting $\mathbf{x}_{k}=\mathbb{B}_{k}\mathbf{y}$
yields
\begin{equation}
\nabla_{k}f(\mathbf{y}_{k},\mathbf{y}_{-k})^{T}(\mathbb{B}_{k}\mathbf{y}-\mathbf{y}_{k})\geq0.\label{eq:cartesian-proof-4}
\end{equation}
Similarly, setting $\mathbf{x}_{k}=\mathbf{y}_{k}$ in (\ref{eq:cartesian-proof-1})
yields
\begin{equation}
\nabla_{k}f(\mathbb{B}_{k}\mathbf{y},\mathbf{y}_{-k})^{T}(\mathbf{y}_{k}-\mathbb{B}_{k}\mathbf{y})\geq0.\label{eq:cartesian-proof-5}
\end{equation}
Adding (\ref{eq:cartesian-proof-4}) and (\ref{eq:cartesian-proof-5}),
we can infer that $(\nabla_{k}f(\mathbf{y})-\nabla_{k}f(\mathbb{B}_{k}\mathbf{y},\mathbf{y}_{-k}))^{T}(\mathbf{y}_{k}-\mathbb{B}_{k}\mathbf{y})\geq0$.
Therefore, we can rewrite (\ref{eq:cartesian-proof-5}) as follows
\begin{align*}
0 & \leq\nabla_{k}f(\mathbb{B}_{k}\mathbf{y},\mathbf{y}_{-k})^{T}(\mathbf{y}_{k}-\mathbb{B}_{k}\mathbf{y})\\
 & =(\nabla_{k}f(\mathbb{B}_{k}\mathbf{y},\mathbf{y}_{-k})-\nabla_{k}f(\mathbf{y})+\nabla_{k}f(\mathbf{y}))^{T}(\mathbf{y}_{k}-\mathbb{B}_{k}\mathbf{y}),
\end{align*}
and thus
\begin{align}
\nabla_{k}f(\mathbf{y})^{T}(\mathbb{B}_{k}\mathbf{y}-\mathbf{y}_{k}) & \leq\nonumber \\
-(\nabla_{k}f(\mathbb{B}_{k}\mathbf{y},\mathbf{y}_{-k})- & \nabla_{k}f(\mathbf{y}))^{T}(\mathbb{B}_{k}\mathbf{y}-\mathbf{y}_{k})\leq0.\label{eq:cartesian-proof-6}
\end{align}
Adding (\ref{eq:cartesian-proof-3}) and (\ref{eq:cartesian-proof-6})
over all $k=1,\ldots,K$ yields
\[
\nabla f(\mathbf{y})^{T}(\mathbb{B}\mathbf{y}-\mathbf{y})={\textstyle \sum_{k=1}^{K}}\nabla_{k}f(\mathbf{y})^{T}(\mathbb{B}_{k}\mathbf{y}-\mathbf{y}_{k})<0.
\]
That is, $\mathbb{B}\mathbf{y}-\mathbf{y}$ is a descent direction
of $f(\mathbf{x})$ at the point $\mathbf{y}$.

The proof of Theorem \ref{thm:convergence} can then be used verbatim
to prove the convergence of the algorithm with the approximate problem
(\ref{eq:jacobi-approximate-problem-cartesian}) and the exact/successive
line search.
\end{IEEEproof}

\end{document}